\documentclass[11pt,reqno]{amsart}
\usepackage{enumerate, latexsym, amsmath, amsfonts, amssymb, amsthm, color}
\def\pmod #1{\ ({\rm{mod}}\ #1)}
\def\Z{\Bbb Z}
\def\N{\Bbb N}

\def\Q{\Bbb Q}

\def\l{\left}
\def\r{\right}
\def\bg{\bigg}
\def\({\bg(}
\def\){\bg)}
\def\t{\text}
\def\f{\frac}

\def\mo{{\rm{mod}\ }}

\def\ord{{\rm ord}}

\def\cs{\ldots}

\def\ls{\leqslant}
\def\gs{\geqslant}
\def\se {\subseteq}
\def\sm{\setminus}
\def\bi{\binom}
\def\al{\alpha}

\def\ve{\varepsilon}

\def\eq{\equiv}

\def\da{\delta}

\def\Proof{\noindent{\it Proof}}
\def\Ack{\medskip\noindent {\bf Acknowledgement}}
\theoremstyle{plain}
\newtheorem{theorem}{Theorem}

\newtheorem{lemma}[theorem]{Lemma}
\newtheorem{corollary}[theorem]{Corollary}
\newtheorem{conjecture}[theorem]{Conjecture}
\theoremstyle{definition}

\theoremstyle{remark}
\newtheorem{remark}[theorem]{Remark}

 \vspace{4mm}

\begin{document}

\hbox{Sci. China Math. 64(2021), no.\,2, 281--306.}
\medskip

\title
[{Further results on Hilbert's Tenth Problem}]
{Further results on Hilbert's Tenth Problem}

\author
[Zhi-Wei Sun] {Zhi-Wei Sun}

\address{Department of Mathematics, Nanjing
University, Nanjing 210093, People's Republic of China}
\email{zwsun@nju.edu.cn}


\begin{abstract} Hilbert's Tenth Problem (HTP) asks for an algorithm to test whether an arbitrary polynomial Diophantine equation
with integer coefficients has solutions over the ring $\Z$ of integers. This was finally solved by Matiyasevich
negatively in 1970. In this paper we obtain some further results on HTP over $\Z$. We prove that there is no algorithm to determine for any $P(z_1,\ldots,z_9)\in\Z[z_1,\ldots,z_9]$
whether the equation $P(z_1,\ldots,z_9)=0$ has {\it integral} solutions with $z_9\gs0$.  Consequently,
there is no algorithm to test whether an arbitrary polynomial Diophantine equation $P(z_1,\ldots,z_{11})=0$ (with integer coefficients) in 11 unknowns has integral solutions,
which provides the best record on the original HTP over $\Z$. We also prove that there is no algorithm to test for any $P(z_1,\ldots,z_{17})\in\Z[z_1,\ldots,z_{17}]$
whether $P(z_1^2,\ldots,z_{17}^2)=0$ has integral solutions, and that there is a polynomial $Q(z_1,\ldots,z_{20})\in\Z[z_1,\ldots,z_{20}]$ such that
$$\{Q(z_1^2,\ldots,z_{20}^2):\ z_1,\ldots,z_{20}\in\Z\}\cap\{0,1,2,\ldots\}$$
coincides with the set of all primes.
\end{abstract}
\maketitle

\section{Introduction}
\setcounter{equation}{0}
\setcounter{theorem}{0}

In 1900, at the Paris conference of International Congress of Mathematicians, D. Hilbert presented 23 famous mathematical problems. He formulated his tenth problem as follows:
\smallskip

 {\it Given a Diophantine equation with any number of unknown quantities and with rational integral numerical coefficients: To devise a process according to which it can be determined in a finite number of operations whether the equation is solvable in rational integers.}
 \smallskip

\noindent In modern language, Hilbert's Tenth Problem (HTP) asks for an effective algorithm to test whether
an arbitrary polynomial equation
$$P(z_1,\cs,z_n)=0$$
(with integer coefficients) has solutions over the ring $\Z$ of integers.
However, at that time the exact meaning of algorithm was not known.

The theory of computability was born
in the 1930s.
A problem or a set is decidable, if and only if
its characteristic function is Turing computable (or recursive).
An {\it r.e. (recursively enumerable) set} is the empty set $\emptyset$ or the range of a recursive function,
and it is also the domain of a partial recursive function.
It is well known that there are nonrecursive r.e. subsets of $\N=\{0,1,2,\ldots\}$.

A relation $R(a_1,\ldots,a_m)$ with $a_1,\ldots,a_m\in\N$ is said to be {\it Diophantine} if
there is a polynomial $P(t_1,\ldots,t_m,x_1,\ldots,x_n)$ with integer coefficients
such that
$$R(a_1,\ldots,a_m)\iff\exists x_1\gs0\ldots\exists x_n\gs0[P(a_1,\ldots,a_m,x_1,\ldots,x_n)=0].$$
(Throughout this paper, variables always range over $\Z$.) A set $\mathcal A\se\N$ is Diophantine if and only if
the predicate $a\in \mathcal A$ is Diophantine. It is easy to see that
any Diophantine set is an r.e. set.

In 1961 Davis et al. \cite{DPR} successfully showed that any r.e. set $\mathcal A$ has an exponential Diophantine representation of the following type:
$$a\in \mathcal A\iff\exists x_1\gs0\ldots\exists x_n\gs0 [P(a,x_1,\ldots,x_n,2^{x_1},\ldots,2^{x_n})=0],$$
where $P$ is a polynomial with integer coefficients.
Recall that the Fibonacci sequence $(F_n)_{n\gs0}$ defined by
$$F_0=0,\ F_1=1,\ \ \t{and}\ F_{n+1}=F_n+F_{n-1}\ (n=1,2,3,\ldots)$$
increases exponentially. In 1970 Yu. Matiyasevich \cite{M70} took the last step to show that
the relation $y=F_{2x}$ (with $x,y\in\N$) is Diophantine. It follows that the exponential relation
$a=b^c$ (with $a,b,c\in\N$, $b>1$ and $c>0$) is Diophantine, i.e. there exists a polynomial $P(a,b,c,x_1,\ldots,x_n)$
with integer coefficients such that
$$a=b^c\iff \exists x_1\gs0\ldots\exists x_n\gs0[P(a,b,c,x_1,\ldots,x_n)=0].$$
This surprising result, together with the important work of Davis et al. \cite{DPR},
leads to the following important result.
\medskip

\noindent {\bf Matiyasevich's Theorem} (See \cite{M70}). {\it Any r.e. set $\mathcal A\se\N$ is Diophantine.}
\smallskip

This famous result is also known as the MDPR theorem named after Matiyasevich, Davis, Putnam and Robinson.

As some r.e. sets are not recursive (cf. \cite[pp.\,140-141]{C80}), Matiyasevich's theorem implies that HTP over $\N$ is undecidable.

Lagrange's four-square theorem in number theory states that any $n\in\N$ can be written as the
sum of four squares. Thus
$P(x_1,\ldots,x_n)=0$ has solutions over $\N$ if and only if
$$P(u_1^2+v_1^2+y_1^2+z_1^2,\ldots,u_n^2+v_n^2+y_n^2+z_n^2)=0$$
has solutions over $\Z$. Now that HTP over $\N$ is undecidable,
so is  HTP over $\Z$ (the original HTP). Thus HTP was finally solved negatively by Matiyasevich in 1970.

It should be mentioned that a whole proof of the unsolvability of HTP is very long
and full of ingenious techniques, see Davis \cite{D73} for a popular introduction,
and Davis, Matiyasevich and Robinson \cite{DMR} for an excellent survey. A machine proof given by Jones
and Matiyasevich \cite{JM} involves clever
arithmetization of register machines (see also Chapter 5 of \cite[pp.\,71-102]{M93}).

For convenience, for a set $S\se\Z$ and a fixed positive integer $n$, we let $\exists^n$ over $S$ denote the set
$$\{\exists x_1\in S\ldots\exists x_n\in S[P(x_1,\ldots,x_n)=0]:\ P(x_1,\ldots,x_n)\in \Z[x_1,\ldots,x_n]\}.$$

Any nonrecursive r.e. set $\mathcal A$ has the following Diophantine representation:
$$a\in \mathcal A\iff\exists x_1\gs0\ldots\exists x_{\nu}\gs0[P(a,x_1,\ldots,x_{\nu})=0],$$
where $P$ is a polynomial with integer coefficients.
Thus $\exists^{\nu}$ over $\N$ is undecidable for some particular number $\nu$. To find the least $\nu$
with $\exists^{\nu}$ over $\N$ undecidable, is a very hard problem.
In the summer of 1970 Matiyasevich announced that $\nu<200$, soon J. Robinson pointed out
that $\nu\ls 35$. Then Matiyasevich and Robinson cooperated in this direction,
in 1973 they \cite{MR} obtained that $\nu\ls 13$, actually they showed that
any diophantine equation over $\N$ can be reduced to one in 13 unknowns.
In 1975, Matiyasevich \cite{M77} announced further that $\nu\ls 9$; a complete proof of this was given by Jones \cite{J82}.

\medskip

\noindent{\bf The 9 Unknowns Theorem} (See \cite{J82}).  {\it $\exists^9$ over $\N$ is undecidable, i.e., there is no algorithm to test whether
$$\exists x_1\gs0\ldots\exists x_9\gs0[P(x_1,\ldots,x_9)=0],$$
where $P(x_1,\ldots,x_9)$ is an arbitrary polynomial in $\Z[x_1,\ldots,x_9]$.}
\medskip

As pointed out by Matiyasevich and Robinson \cite{MR}, if $a_0,a_1,\ldots,a_n$ and $z$ are integers
with $a_0z\not=0$ and $\sum_{i=0}^na_iz^{n-i}=0$, then
$$|z|^n\ls|a_0z^n|\ls\sum_{i=1}^n|a_i|\cdot|z|^{n-i}\ls \sum_{i=1}^n|a_i|\cdot|z|^{n-1}\ \t{and hence}\ |z|\ls \sum_{i=1}^n|a_i|.$$
Thus $\exists$ over $\N$ and $\exists$ over $\Z$ are decidable.
It is not known whether $\exists^2$ over $\N$ or $\exists^2$ over $\Z$ is decidable, though A. Baker \cite{B68} showed that if $F[x,y]\in\Z[x,y]$ is irreducible, homogeneous and of degree at least three then for any $m\in\Z$ there is an effective algorithm to find integral solutions of the equation $F(x,y)=m$.
Baker \cite{B68}, Matiyasevich and Robinson \cite{MR} believed that $\exists^3$ over $\N$ is undecidable.

As the original HTP is about integral solutions of polynomial Diophantine equations, it is natural to ask for the smallest $\mu\in\Z^+=\{1,2,3,\ldots\}$
such that $\exists^{\mu}$ over $\Z$ is undecidable.
In view of Lagrange's four-square theorem,
if $\exists^n$ over $\N$ is undecidable, then so is
$\exists^{4n}$ over $\Z$. This can be made better.
By the Gauss-Legendre theorem on sums of three squares (cf. \cite[pp.\, 17-23]{N96}), the number $4m+1$ with $m\in\N$ can be written as the sum of
two even squares and an odd square. It follows that for any integer $m$ we have
\begin{equation}\label{1.1}m\gs0\iff \exists x\exists y\exists z[m=x^2+y^2+z^2+z].\end{equation}
Therefore the undecidability of $\exists^n$ over $\N$ implies
the undecidability of $\exists^{3n}$ over $\Z$, thus Tung \cite{T85}
obtained the undecidability of $\exists^{27}$ over $\Z$ from the
9 unknowns theorem. Tung \cite{T85} asked whether 27 here can be replaced by a smaller number.
In 1992 Sun \cite{S92b} showed for any $n\in\Z^+$ that if $\exists^n$ over $\N$ is undecidable then so is $\exists^{2n+2}$ over $\Z$,
and thus he obtained the undecidability of $\exists^{20}$ over $\Z$ from the 9 unknowns theorem.
The author announced in \cite{S92a, S92b} that $\exists^{11}$ over $\Z$ is undecidable,
however the whole sophisticated proof appeared in his PhD thesis
 \cite{S-PhD} has not been published before.

HTP over the field $\Q$ of rational numbers remains open, but Robinson \cite{Ro49} showed that
the first-order theory of rational numbers is undecidable (see D. Flath and S. Wagon \cite{FW} for an excellent introduction and J. Koenigsmann \cite{K} for recent progress). There are also lots of research works on extended HTP over various rings and fields (see, e.g., J. Denef \cite{De75, De78}, Denef and L. Lipshitz \cite{DL}, and A. Shlapentokh \cite{Sh}).

For the extended HTP over a ring $R$ containing $\Z$,
the usual strategy to obtain its undecidability is as follows: Prove that $\Z$ is Diophantine over $R$ and then use the result that HTP over $\Z$ is undecidable.
Thus, to find a small positive integer $k$ with $\exists^k$ over $R$ undecidable, depends heavily on
the undecidability of $\exists^{\mu}$ over $\Z$ ({\it not} $\N$) with $\mu$ as small as possible. In this sense, to find a small number
$\mu$ with $\exists^{\mu}$ over $\Z$ undecidable is quite important and very useful.

In this paper we focus on HTP over $\Z$. Now we state our first theorem which implies the undecidability of $\exists^{11}$ over $\Z$.

\begin{theorem}\label{Th1.1} Let $\mathcal A\se\N$ be any r.e. set.

{\rm (i)} There is a polynomial $P_{\mathcal A}(z_0,z_1,\ldots,z_9)$
with integer coefficients such that for any $a\in\N$ we have
\begin{equation}\label{1.2}\exists z_1\ldots\exists z_8\exists z_9\gs0[P_{\mathcal A}(a,z_1,\ldots,z_9)=0]
\Longrightarrow a\in \mathcal A,\end{equation}
and
\begin{equation}\label{1.3}a\in\mathcal A\Longrightarrow \forall Z>0\exists z_1\gs Z\ldots\exists z_9\gs Z [P_{\mathcal A}(a,z_1,\ldots,z_9)=0].\end{equation}

{\rm (ii)} There is a polynomial $Q_{\mathcal A}(z_0,z_1,\ldots,z_{10})$ with integer coefficients
such that for any $a\in\N$ we have
\begin{equation}\label{1.4}a\in\mathcal A\iff\exists z_1\ldots\exists z_9\exists z_{10}\not=0[Q_{\mathcal A}(a,z_1,\ldots,z_{10})=0].\end{equation}
\end{theorem}

Since there are nonrecursive r.e. sets, with the aid of (\ref{1.1}) and Tung's observation (see \cite{T85}) that
\begin{equation}\label{1.5}m\in\Z\sm\{0\}\iff\exists x\exists y [m=(2x+1)(3y+1)],\end{equation}
we immediately get the following consequence.

\begin{corollary}\label{Cor1.1} {\rm (i) (The 11 Unknowns Theorem)} $\exists^{11}$ over $\Z$ is undecidable. Moreover, there is no algorithm to determine whether the equation
\begin{equation}\label{1.6}P(z_1,\ldots,z_8,z_9^2+z_{10}^2+z_{11}^2+z_{11})=0\end{equation}
 has solutions over $\Z$ for an arbitrary polynomial $P(z_1,\ldots,z_9)\in\Z[z_1,\ldots,z_9]$. Also, there is no algorithm to determine whether the equation
\begin{equation}\label{1.7}Q(z_1,\ldots,z_9,(2z_{10}+1)(3z_{11}+1))=0\end{equation}
 has integral solutions for an arbitrary polynomial $Q(z_1,\ldots,z_{10})\in\Z[z_1,\ldots,z_{10}]$.

{\rm (ii)} There is no algorithm to test whether
\begin{equation}\label{1.8}\forall Z>0\exists z_1\gs Z\ldots\exists z_9\gs Z[P(z_1,\ldots,z_{9})=0],\end{equation}
where $P(z_1,\ldots,z_9)$ is an arbitrary polynomial in $\Z[z_1,\ldots,z_9]$.
\end{corollary}
\begin{remark}\label{Rem1.2} Corollary \ref{Cor1.1}(i) provides the best record on the original HTP over $\Z$.
In the author's opinion, this can hardly be improved in a near future.
\end{remark}

In number theory, a subset $S$ of $\N$ is called an asymptotic additive base of order $h$ if all sufficiently large integers can be written as $a_1+\ldots+a_h$
with $a_1,\ldots,a_h\in S$. From Theorem \ref{Th1.1}(i) we see that if $S\se\N$ is an asymptotic additive base of order $h$
then $\exists^{9h}$ over $S$ is undecidable. Thus, $\exists^{9G(k)}$ over $\{m^k:\ m\in\N\}$ is undecidable for every $k=2,3,4,\ldots$, where
$G(k)$ associated with Waring's problem denotes the least positive integer $s$ such that any sufficiently large integer can be written as $x_1^k+\ldots+x_s^k$
with $x_1,\ldots,x_s\in\N$. It is known that $G(2)=4$, $G(3)\ls 7$ and $G(4)=16$ (cf. \cite{V97}).
\medskip

\begin{corollary}\label{Cor1.2} {\rm (i)} $\forall^{9}\exists^3$ over $\Z$ is undecidable, i.e., there is no algorithm to test whether
\begin{equation}\label{1.9}\forall z_1\ldots\forall z_{9}\exists x\exists y\exists z[P(z_1,\ldots,z_{9},x,y,z)=0],\end{equation}
where $P$ is an arbitrary polynomial of $12$ variables with integer coefficients.

 {\rm (ii)} $\forall^{10}\exists^2$ over $\Z$ is undecidable, i.e., there is no algorithm to test whether
\begin{equation}\label{1.10}\forall z_1\ldots\forall z_{10}\exists x\exists y [Q(z_1,\ldots,z_{10},x,y)=0],\end{equation}
where $Q$ is an arbitrary polynomial of $12$ variables with integer coefficients.
\end{corollary}
\begin{remark}\label{Rem1.3} In 1981 Jones \cite{J81} obtained the decidability of $\forall \exists$ over $\N$ as well as some other undecidable results over $\N$.
In 1987 Tung \cite{T87} proved for each $n\in\Z^+$ that $\forall^n\exists $ over $\Z$ is co-NP-complete. Tung \cite{T87} also showed that
$\forall^{27}\exists^2 $ over $\Z$ is undecidable, and asked whether $27$ here can be replaced by a smaller number.
\end{remark}

Our next theorem is related to polygonal numbers.
Recall that triangular numbers have the form $T_x=x(x+1)/2$ with $x\in\Z$, generalized pentagonal numbers are those integers $p_5(x)=x(3x-1)/2$ with $x\in\Z$,
and generalized octagonal numbers are those $p_8(x)=x(3x-2)$ with $x\in\Z$. Polygonal numbers of order four coincide with squares of integers.

\begin{theorem}\label{Th1.2} Let $\mathcal A$ be any r.e. subset of $\N$.
Then there is a polynomial $P_4(z_0,z_1,\ldots,z_{17})$
with integer coefficients such that for any $a\in\N$ we have
\begin{equation}\label{1.11}a\in\mathcal A\iff \exists z_1\in\square\ldots\exists z_{17}\in\square[P_4(a,z_1,\ldots,z_{17})=0],\end{equation}
where $\square$ denotes the set of all integer squares.
Also, there are polynomials
$$P_3(z_0,z_1,\ldots,z_{18}),\ P_5(z_0,z_1,\ldots,z_{18}),\ P_8(z_0,z_1,\ldots,z_{18})$$
with integer coefficients such that for any $a\in\N$ we have
\begin{equation}\label{1.12}\begin{aligned} a\in\mathcal A\iff& \exists z_1\in{\rm Tri}\ldots\exists z_{18}\in{\rm Tri}[P_3(a,z_1,\ldots,z_{18})=0]
\\\iff&\exists z_1\in{\rm Pen}\ldots\exists z_{18}\in{\rm Pen}[P_5(a,z_1,\ldots,z_{18})=0]
\\\iff&\exists z_1\in{\rm Octa}\ldots\exists z_{18}\in{\rm Octa}[P_8(a,z_1,\ldots,z_{18})=0],
\end{aligned}\end{equation}
where
 \begin{gather*}{\rm Tri}=\{T_x:\ x\in\Z\}=\{x(2x+1):\ x\in\Z\},
  \\ {\rm Pen}=\{p_5(x):\ x\in\Z\}\ \ \t{and}\ \ {\rm Octa}=\{p_8(x):\ x\in\Z\}.
  \end{gather*}
\end{theorem}

Clearly Theorem \ref{Th1.2} has the following consequence.

\begin{corollary}\label{Cor1.3} $\exists^{17}$ over $\square$,
$\exists^{18}$ over ${\rm Tri}$, $\exists^{18}$ over ${\rm Pen}$,
and $\exists^{18}$ over ${\rm Octa}$ are all undecidable.
\end{corollary}

Motivated by Corollary \ref{Cor1.3}, we formulate the following conjecture.
\begin{conjecture}\label{Con1.1} $\exists^3$ over $\square$ is undecidable, i.e., there is no algorithm to determine for any $P(x,y,z)\in\Z[x,y,z]$
whether the equation $P(x^2,y^2,z^2)=0$ has integral solutions.
\end{conjecture}

Using Theorems \ref{Th1.1} and \ref{Th1.2}, we deduce the following result.

\begin{theorem}\label{Th1.3} {\rm (i)} Let $\mathcal A\se\N$ be any r.e. set. Then there is a polynomial $P(z_1,\ldots,z_{14})$ with integer coefficients such that
\begin{equation}\label{1.13}\mathcal A=\N\cap\{P(z_1,\ldots,z_{14}):\ z_1,\ldots,z_{14}\in\Z\}.\end{equation}
Also, there are polynomials
$$Q_4(z_1,\ldots,z_{21}),\ Q_3(z_1,\ldots,z_{21}),\ Q_5(z_1,\ldots,z_{21}),\ Q_8(z_1,\ldots,z_{22})$$
with integer coefficients such that
\begin{equation}\label{1.14}\begin{aligned} \mathcal A=&\N\cap\{Q_4(z_1,\ldots,z_{21}):\ z_1,\ldots,z_{21}\in\square\}
\\=&\N\cap\{Q_3(z_1,\ldots,z_{21}):\ z_1,\ldots,z_{21}\in{\rm Tri}\}
\\=&\N\cap\{Q_5(z_1,\ldots,z_{21}):\ z_1,\ldots,z_{21}\in{\rm Pen}\}
\\=&\N\cap\{Q_8(z_1,\ldots,z_{22}):\ z_1,\ldots,z_{22}\in{\rm Octa}\}.
\end{aligned}\end{equation}

{\rm (ii)} Let $\mathcal P$ be the set of all primes. There are polynomials $\hat P(z_1,\ldots,z_{20})$ and $\tilde P(z_1,\ldots,z_{21})$
with integer coefficients such that
\begin{equation}\label{1.15}\mathcal P=\N\cap\{\hat P(z_1^2,\ldots,z_{20}^2):\ z_1,\ldots,z_{20}\in\Z\}\end{equation}
and
\begin{equation}\label{1.16}\mathcal P=\N\cap\{\tilde P(z_1(3z_1+2),\ldots,z_{21}(3z_{21}+2)):\ z_1,\ldots,z_{21}\in\Z\}.\end{equation}
\end{theorem}
\begin{remark}\label{Rem1.4} Matiyasevich \cite{M81} constructed a polynomial
$P(x_1,\ldots,x_{10})$ with integer coefficients such that
$$\mathcal P=\N\cap\{P(x_1,\ldots,x_{10}):\ x_1,\ldots,x_{10}\in\N\}.$$
\end{remark}

To give detailed proofs of Theorems \ref{Th1.1} and \ref{Th1.3}, we utilize some basic ideas in Matiyasevich's proof of the 9 unknowns theorem (cf. \cite{J82})
as well as the earlier coding idea of Matiyasevich and Robinson \cite{MR} on reduction of unknowns, and we also
overcome various new technical difficulties caused by avoiding natural number variables, and employ some recent results of the author on polygonal numbers.
Our starting point is the use of Lucas sequences with integer indices.

Let $A$ and $B$ be integers. The usual Lucas sequence $u_n=u_n(A,B)\ (n=0,1,2,\ldots)$ and its companion
$v_n=v_n(A,B)\ (n=0,1,2,\ldots)$ are defined as follows:
$$u_0=0,\ u_1=1,\ \t{and}\ u_{n+1}=Au_n-Bu_{n-1}\ (n=1,2,3,\ldots);$$
and
$$v_0=2,\ v_1=A,\ \t{and}\ v_{n+1}=Av_n-Bv_{n-1}\ (n=1,2,3,\ldots).$$
Note that $u_n(2,1)=n$, $u_n(1,-1)=F_n$ and $u_n(3,1)=F_{2n}$ for all $n\in\N$.
Let
$$\al=\f{A+\sqrt{\Delta}}2\ \t{and}\ \ \beta=\f{A-\sqrt{\Delta}}2$$
be the two roots of the quadratic equation $x^2-Ax+B=0$ where $\Delta=A^2-4B$.
It is well known that
\begin{equation}\label{1.17}(\al-\beta)u_n=\al^n-\beta^n, \ v_n=\al^n+\beta^n\ \t{and}\ v_n^2-\Delta u_n^2=4B^n\end{equation}
for all $n\in\N$ (see, e.g., [Ri89, pp.\,41-42]). If $u_n\gs0$ for all $n\in\N$, then $A=u_2\gs0$ and $\Delta\gs0$ (otherwise $u_{n+1}^2-u_nu_{n+2}=B^n>0$
and the decreasing sequence $(u_{n+1}/u_n)_{n\gs1}$ has a limit which should be a real root of the equation $x^2-Ax+B=0$).
Conversely, if $A\gs0$ and $\Delta\gs0$ then $u_n\gs0$ for all $n\in\N$, which can be easily shown. When $\Delta\gs0$,
the sequence $(u_n)_{n\gs0}$ is strictly increasing if and only if $A>1$ (cf. \cite[Lemma 4]{S92a}).

We actually only need Lucas sequences with $B=1$. In this case, we extend the sequences $u_n=u_n(A,1)$ and $v_n=v_n(A,1)$ to integer indices by letting
\begin{equation}\label{1.18}u_0=0,\ u_1=1,\ \t{and}\ u_{n-1}+u_{n+1}=Au_n\ \t{for all}\ n\in\Z,\end{equation}
and
\begin{equation}\label{1.19}v_0=2,\ v_1=A,\ \t{and}\ v_{n-1}+v_{n+1}=Av_n\ \t{for all}\ n\in\Z.\end{equation}
It is easy to see that
\begin{equation}\label{1.20}u_{-n}(A,1)=-u_n(A,1)=(-1)^nu_n(-A,1)\end{equation}
and $v_{-n}(A,1)=v_n(A,1)=(-1)^nv_n(-A,1)$ for all $n\in\Z$.
For the relation $C=u_B(A,1)$ with $A,B,C\in\Z$, the author studied its Diophantine representations over $\Z$ in the published paper \cite{S92a}.
This laid the initial foundation for our work in this paper.

We provide some lemmas on $p$-adic expansions in the next section and then show an auxiliary theorem in Section 3.
In Section 4 we work with Lucas sequences and prove two key theorems on Diophantine representations.
In Section 5 we prove Theorem \ref{Th1.1} and Corollary \ref{Cor1.2}. Section 6 is devoted to our proofs of Theorems \ref{Th1.2} and \ref{Th1.3}.

Throughout this paper, we adopt the notation
$$ p\uparrow:=\{p^n:\ n\in\N\}\ \quad \t{for}\ p\in\Z^+.$$
For $c,d\in\Z$ we define $[c,d):=\{m\in\Z:\ c\ls m<d\}$. For a prime $p$ and a nonzero integer $m$, we use $\ord_p(m)$
to denote the $p$-adic order of $m$ at $p$, i.e., the largest $a\in\N$ with $p^a$ dividing $m$.
All the 26 capital Latin letters
$A,B,\ldots,Y,Z$ will be used in our proofs of Theorems \ref{Th1.1} and \ref{Th1.2}, and each of them has a special meaning.

\section{Some lemmas on $p$-adic expansions}
\setcounter{equation}{0}
\setcounter{theorem}{0}

Let $p>1$ be an integer. Any $n\in\N$ has a unique $p$-adic expansion
$$\sum_{i=0}^\infty a_ip^i\ \ \t{with}\ \ a_i\in[0,p)=\{0,1,\ldots,p-1\},$$
where $a_j=0$ for all sufficiently large values of $j$.
Let $$\sigma_p(n):=\sum_{i=0}^\infty a_i$$
 be the sum of all digits in the $p$-adic (or base $p$) expansion of $n$.
 Since $a_i=\lfloor n/p^i\rfloor-p\lfloor a/p^{i+1}\rfloor$, we see that
 \begin{equation}\label{2.1}\sigma_p(n)=\sum_{i=0}^\infty\l(\l\lfloor\f{n}{p^i}\r\rfloor-p\l\lfloor\f n{p^{i+1}}\r\rfloor\r)=n-(p-1)\sum_{i=1}^\infty\l\lfloor\f n{p^i}\r\rfloor\end{equation}
 as first observed by Legendre (cf. \cite[p.\,22]{Ri89}).

If $p$ is a prime, then
$$\ord_p(n!)=\sum_{i=1}^\infty\l\lfloor\f n{p^i}\r\rfloor\ \quad\t{for all}\ n\in\N.$$
Combining this well-known result with (\ref{2.1}), we immediately get the following result essentially due to Kummer (cf. \cite[pp.\,23-24]{Ri89})

\begin{lemma}\label{Lem2.1} Let $a,b\in\N$ and let $p$ be a prime.
Let $\tau_p(a,b)$ denote the number of carries occurring in the addition of $a$ and $b$ in base $p$. Then
\begin{equation}\label{2.2}\tau_p(a,b)=\ord_p\bi{a+b}a=\f{\sigma_p(a)+\sigma_p(b)-\sigma_p(a+b)}{p-1}.
\end{equation}
\end{lemma}

With the aid of Lemma \ref{Lem2.1}, we deduce the following lemma.

\begin{lemma}\label{Lem2.2} Let $p$ be a prime, and let $P\in p\uparrow$, $N\in P\uparrow$ and $S,T\in[0,N)$. Then
\begin{equation}\label{2.3}\tau_p(S,T)=0\iff N^2\mid\bi {P\f{N-1}{P-1}R}{\f{N-1}{P-1}R},\end{equation}
where $R:=(S+T+1)N+T+1$.
\end{lemma}
\Proof. Write $N=p^n$ with $n\in\N$. By Lemma \ref{Lem2.1}, we have
\begin{align*}&N^2\mid\bi{P\f{N-1}{P-1}R}{\f{N-1}{P-1}R}
\\\iff& \tau_p\l((N-1)R,\f{N-1}{P-1}R\r)\gs 2n
\\\iff&\sigma_p((N-1)R)+\sigma_p\l(\f{N-1}{P-1}R\r)-\sigma_p\l(P\f{N-1}{P-1}R\r)\gs 2n(p-1).
\end{align*}
Clearly, $\sigma_p(Pm)=\sigma_p(m)$ for any $m\in\N$. Note that
$$(N-1)R=(S+T)N^2+(N-1-S)N+N-1-T.$$
Thus
$$\sigma_p((N-1)R)=\sigma_p(S+T)+\sigma_p(N-1-S)+\sigma_p(N-1-T).$$
As $N-1=\sum_{0\ls i<n}(p-1)p^i$, we see that
$$\sigma_p(N-1-S)=n(p-1)-\sigma_p(S)\ \ \t{and}\ \ \sigma_p(N-1-T)=n(p-1)-\sigma_p(T).$$
Therefore
\begin{align*}&N^2\mid\bi{P\f{N-1}{P-1}R}{\f{N-1}{P-1}R}
\\\iff&\sigma_p((N-1)R)\gs2n(p-1)
\\\iff&\sigma_p(S+T)+(n(p-1)-\sigma_p(S))+(n(p-1)-\sigma_p(T))\gs 2n(p-1)
\\\iff& \sigma_p(S)+\sigma_p(T)\ls \sigma_p(S+T).
\end{align*}
By Lemma \ref{Lem2.1},
$$\sigma_p(S)+\sigma_p(T)\ls \sigma_p(S+T)\iff \tau_p(S,T)\ls0\iff \tau_p(S,T)=0.$$
So the desired result follows. \qed
\medskip

\begin{remark}\label{Rem2.1} Lemma \ref{Lem2.2} in the case $P=p=2$ appeared in \cite[Lemma 2.16]{J82}.
\end{remark}

\begin{lemma}\label{Lem2.3} Let $p>1$ be an integer and let $b,B\in p\uparrow$ with $b\ls B$.
Let $n_1,\ldots,n_k\in\N$ with $n_1<\ldots<n_k$. Suppose that $C\in\Z^+$ with $b\ls C/B^{n_k}\ls B$.
Then
$$c=\sum_{i=1}^kz_iB^{n_i}\ \t{for some}\ z_1,\ldots,z_k\in[0,b)\iff c\in[0,C)\land \tau_p(c,M)=0,$$
where $M=\sum_{j=0}^{n_k}m_jB^j$, and
$$m_j=\begin{cases} B-b&\t{if}\ j\in\{n_s:\ s=1,\ldots,k\},\\B-1&\t{otherwise}.\end{cases}$$
\end{lemma}
\Proof. If $c=\sum_{i=1}^kz_iB^{n_i}$ for some $z_1,\ldots,z_k\in[0,b)$, then
$$0\ls c\ls\sum_{i=1}^k(b-1)B^{n_i}\ls(b-1)B^{n_k}+\sum_{j=0}^{n_k-1}(B-1)B^j=bB^{n_k}-1<C.$$

Let $c\in[0,C)$. As $c\ls B^{n_k+1}-1=\sum_{j=0}^{n_k}(B-1)B^j$, we can write
$c=\sum_{j=0}^{n_k}c_jB^j$ with $c_j\in[0,B)$. If $B=1$ then $c=0$ and $\tau_p(c,M)=0$.
Since $b,B\in p\uparrow$, when $B>1$ we have
\begin{align*}\tau_p(c,M)=0\iff& \tau_p(c_j,m_j)=0\ \t {for all}\ j=0,\ldots,n_k
\\\iff&c_{n_i}<b\ \t{for all}\ i=1,\ldots,k,\ \t{and}\ c_j=0\ \t{for}\ j\not\in\{n_s:\ 1\ls s\ls k\}
\\\iff&c=\sum_{i=1}^kz_iB^{n_i}\ \t{for some}\ z_1,\ldots,z_k\in[0,b).
\end{align*}
This concludes the proof. \qed
\medskip

\begin{remark}\label{Rem2.2} Lemma \ref{Lem2.3} and the following Lemma \ref{Lem2.4} utilize some coding ideas of Matiyasevich and Robinson (cf. \cite[Section 6]{MR}
and \cite[Section 3]{J82}) who worked in the case $p=2$.
\end{remark}

\begin{lemma}\label{Lem2.4} Let $\da\in\Z^+$, $z_0,\ldots,z_{\nu}\in\N$ and
$$P(z_0,z_1,\ldots,z_{\nu})=\sum_{i_0,\ldots,i_{\nu}\in\N\atop i_0+\ldots+i_{\nu}\ls\da} a_{i_0,\ldots,i_{\nu}}z_0^{i_0}\cdots z_{\nu}^{i_{\nu}}$$
with $a_{i_0,\ldots,i_{\nu}}\in\Z$ and $|a_{i_0,\ldots,i_{\nu}}|\ls L\in\Z^+$.
Let $p>1$ be an integer, and let $B,X\in p\uparrow$ with
$$B>X>\da!L(1+z_0+z_1+\cdots+z_{\nu})^\da.$$
Let $n_i=(\da+1)^i$ for $i=0,1,2,\ldots$. Set $c=1+\sum_{i=0}^{\nu}z_iB^{n_i}$ and
\begin{align*}K=&c^\da\sum_{i_0,\ldots,i_{\nu}\in\N\atop i_0+\cdots+i_{\nu}\ls\da} i_0!\cdots i_{\nu}!(\da-i_0-\cdots-i_{\nu})!a_{i_0,\ldots,i_{\nu}}B^{n_{\nu+1}-\sum_{s=0}^{\nu}i_s n_s}
\\&+X\sum_{i=0}^{(2\da+1)n_{\nu}}B^i.
\end{align*}
Then $B^{(2\da+1)n_{\nu}}<K<B^{(2\da+1)n_{\nu}+1}$, and
\begin{equation}\label{2.4}P(z_0,\ldots,z_{\nu})=0\iff \tau_p(K,(X-1)B^{n_{\nu+1}})=0.\end{equation}
\end{lemma}
\Proof. Write
$$C(x):=\(1+\sum_{i=0}^{\nu}z_ix^{n_i}\)^\da=\sum_{i=0}^{\da n_{\nu}}c_ix^i$$
and
$$D(x):=\sum_{i_0,\ldots,i_{\nu}\in\N\atop i_0+\cdots+i_{\nu}\ls\da} i_0!\cdots i_{\nu}!(\da-i_0-\cdots-i_{\nu})!a_{i_0,\ldots,i_{\nu}}x^{n_{\nu+1}-\sum_{s=0}^\nu i_sn_s}
=\sum_{j=0}^{n_{\nu+1}}d_jx^j.$$
Clearly, $C(B)=c^\da$ and
$$K=C(B)D(B)+X\sum_{i=0}^{(2\da+1)n_{\nu}}B^i=\sum_{k=0}^{(2\da+1)n_{\nu}}e_kB^k$$
with
$$e_k=X+\sum_{0\ls i\ls\da n_{\nu}\atop {0\ls j\ls n_{\nu+1}\atop i+j=k}} c_id_j.$$
For $i_0,\ldots,i_{\nu}\in\N$ with $i_0+\ldots+i_{\nu}\ls\da$, the multi-nomial coefficient
$$\bi{\da}{i_0,\ldots,i_{\nu},\da-i_0-\cdots-i_{\nu}}=\f{\da!}{i_0!\cdots i_{\nu}!(\da-i_0-\cdots-i_{\nu})!}$$
is a positive integer and hence
\begin{equation}\label{2.5}i_0!\cdots i_{\nu}!(\da-i_0-\cdots-i_{\nu})!\ls\da!.\end{equation}
As $|d_j|\ls\da!L$ for all $j=0,\ldots,n_{\nu+1}$, we have
$$|e_k-X|\ls\da!L\sum_{i=0}^{\da n_{\nu}}c_i=\da!LC(1)=\da!L(1+z_0+\cdots+z_{\nu})^\da<X$$
and hence $0<e_k<2X\ls pX\ls B$. It follows that
\begin{align*}B^{(2\da+1)n_{\nu}}<\sum_{k=0}^{(2\da+1)n_{\nu}}B^k\ls &K=\sum_{k=0}^{(2\da+1)n_{\nu}}e_kB^k
\ls (B-1)\sum_{k=0}^{(2\da+1)n_{\nu}}B^k<B^{(2\da+1)n_{\nu}+1}.
\end{align*}

 By the multi-nomial theorem,
$$C(x)=\sum_{i_0,\ldots,i_{\nu}\in\N\atop i_0+\cdots+i_{\nu}\ls\da}\f{\da!}{i_0!\cdots i_{\nu}!(\da-i_0-\cdots-i_{\nu})!}z_0^{i_0}\cdots z_{\nu}^{i_{\nu}}x^{\sum_{s=0}^\nu i_sn_s}.$$
Recall that $n_s=(\da+1)^s$. The coefficient of $x^{n_{\nu+1}}$ in the expansion of $C(x)D(x)$ coincides with
$$\sum_{i_0,\ldots,i_{\nu}\in\N\atop i_0+\cdots+i_{\nu}\ls\da}\da!z_0^{i_0}\cdots z_{\nu}^{i_{\nu}}a_{i_0,\ldots,i_{\nu}}=\da!P(z_0,\ldots,z_{\nu}),$$
and hence
$$-X<\da!P(z_0,\ldots,z_{\nu})=e_{n_{\nu+1}}-X<X.$$
As $K=\sum_{k=0}^{(2\da+1)n_{\nu}}e_kB^k$ with $0<e_k<B$, and $0\ls X-1<X<B$, we have
$$\tau_p(K,(X-1)B^{n_{\nu+1}})=0\iff \tau_p(e_{n_{\nu+1}},X-1)=0$$
since $B\in p\uparrow$. As $X\in p\uparrow$ and $e_{n_{\nu+1}}\in[1,2X)$, we see that
$$\tau_p(e_{n_{\nu+1}},X-1)=0\iff e_{n_{\nu+1}}=X\iff P(z_0,\ldots,z_{\nu})=0.$$
Therefore (\ref{2.4}) does hold. \qed

\section{An auxiliary theorem}
\setcounter{equation}{0}
\setcounter{theorem}{0}

In this section, we employ lemmas in Section 2 to establish the following auxiliary result which is indispensable for our later proofs of Theorems \ref{Th1.1} and \ref{Th1.2}.

\begin{theorem}\label{Th3.1} Let $\mathcal A\se\N$ be a Diophantine set, and let $p$ be a prime. Then, for each $a\in\N$, we have
\begin{equation}\label{3.1}
a\in\mathcal A\Rightarrow\forall Z>0\exists f\gs Z\exists g\in [b,\mathcal C)\(b\in\square\land b\in p\uparrow\land Y\mid\bi{pX}X\)\end{equation}
and
\begin{equation}\label{3.2}
\exists f\not=0\exists g\in[0,2\mathcal C)\l(b\in\square\land b\in p\uparrow\land Y\mid\bi{pX}X\r)\Rightarrow a\in \mathcal A,\end{equation}
where
\begin{equation}\label{3.3}b:=1+(p^2-1)(ap+1)f,\end{equation}
$\mathcal C=p^{\al_1p}b^{\al_2}$ for some $\al_1,\al_2\in\Z^+$ only depending on $\mathcal A$, and $X$ and $Y$ are suitable polynomials in $\Z[a,f,g]$
such that
if $a\in\N$, $f\in\Z\sm\{0\}$, $b\in \square$ and $0\ls g<2\mathcal C$ then
\begin{equation}\label{3.4}p+1\mid X,\ \ X\gs3b\ \ \t{and}\ \ Y\gs \max\{b,p^{4p}\}.\end{equation}
\end{theorem}
\Proof. As the set $\mathcal A$ is Diophantine, there is a polynomial
$P(z_0,z_1,\ldots,z_{\nu})$ with integer coefficients such that for any $a\in\N$ we have
$$a\in \mathcal A\iff\exists z_1\gs0\ldots\exists z_{\nu}\gs0[P(a,z_1,\ldots,z_{\nu})=0].$$
Thus
$$a\in \mathcal A\iff \exists z_1\gs0\ldots\exists z_{\nu+1}\gs0[\bar P(a,z_1,\ldots,z_{\nu+1})=0],$$
where
$$\bar P(z_0,z_1,\ldots,z_{\nu+1})=P(z_0,z_1,\ldots,z_{\nu})^2+(z_{\nu+1}-1)^2$$
with $\bar P(a,0,\ldots,0)=P(a,0,\ldots,0)^2+(0-1)^2>0.$ Without loss of generality, we simply assume that
$P(a,0,\ldots,0)>0$ for all $a\in\N$. Write
$$P(z_0,\ldots,z_{\nu})=\sum_{i_0,\ldots,i_{\nu}\in\N\atop i_0+\ldots+i_{\nu}\ls\da} a_{i_0,\ldots,i_{\nu}}z_0^{i_0}\ldots z_{\nu}^{i_{\nu}}$$
with $a_{i_0,\ldots,i_{\nu}}\in\Z$, where $\da\in\Z^+$. For
$$\mathcal L:=\max_{i_0,\ldots,i_{\nu}\in\N\atop i_0+\ldots+i_{\nu}\ls\da}|a_{i_0,\ldots,i_{\nu}}|,$$
we obviously have
$$\mathcal L\gs a_{0,\ldots,0}=P(0,\ldots,0)>0.$$

Let $a\in\N$. As $p$ is relatively prime to $(p^2-1)(ap+1)$, by Euler's theorem we have
$$p^{\varphi((p^2-1)(ap+1))}\eq1\pmod {(p^2-1)(ap+1)},$$  where $\varphi$ is Euler's totient function.
Let $Z\in\Z^+$. If $a\in\mathcal A$, then $P(a,z_1,\ldots,z_{\nu})$ $=0$ for some $z_1,\ldots,z_{\nu}\in\N$, hence we may take a sufficiently large integer $n>0$
such that
$$b_0:=p^{2n\varphi((p^2-1)(ap+1))}>\max\{z_1,\ldots,z_{\nu},1+(p^2-1)(ap+1)Z\},$$
and this $b_0$ is a square and it can be written as $1+(p^2-1)(ap+1)f_0$ with $f_0\in\Z$ and $f_0\gs Z$.

Now fix $a\in\N$, and suppose that $f\in\Z\sm\{0\}$ and $b=1+(p^2-1)(ap+1)f\in\square$. Clearly, $f>0$ and hence $b\gs ap+1>a$.
Note that $0<c:=(\nu+1)b<(\nu+2)b-a$.
Take a positive integer $\al$ with
$$\beta:=p^{\al p}>(\nu+2)^{\da}\da!p\mathcal L.$$
Then
$$\mathcal B:=\beta b^\da>(\nu+2)^\da\da!p\mathcal Lb^\da\gs(a+c+1)^\da\da!p\mathcal L\gs p+(a+c)^\da \da!p\mathcal L$$
and
\begin{equation}\label{3.5}\f{\mathcal B}p>(a+c)^\da \da!\mathcal L\gs\da!\mathcal L(1+a+\nu(b-1))^\da.
\end{equation}
Define
$$D(\mathcal B):=\sum_{i_0,\ldots,i_{\nu}\in\N\atop i_0+\cdots+i_{\nu}\ls\da} i_0!\cdots i_{\nu}!(\da-i_0-\cdots-i_{\nu})!a_{i_0,\ldots,i_{\nu}}
\mathcal B^{(\da+1)^{\nu+1}-\sum_{s=0}^\nu i_s(\da+1)^s}.$$
In view of (2.5), we have
\begin{align*}&\left|D(\mathcal B)-\da!a_{0,\ldots,0}\mathcal B^{(\da+1)^{\nu+1}}\right|
\\\ls&\sum_{i_0,\ldots,i_{\nu}\in\N\atop 0<i_0+\cdots+i_{\nu}\ls\da} \da!|a_{i_0,\ldots,i_{\nu}}|\mathcal B^{(\da+1)^{\nu+1}-\sum_{s=0}^\nu i_s(\da+1)^s}
\\\ls&\da!\mathcal L\sum_{r=0}^{(\da+1)^{\nu+1}-1}\mathcal B^r\ls(\mathcal B-1)\sum_{r=0}^{(\da+1)^{\nu+1}-1}\mathcal B^r<B^{(\da+1)^{\nu+1}}
\end{align*}
and hence
\begin{equation}\label{3.6}
D(\mathcal B)>(\da!a_{0,\ldots,0}-1)\mathcal B^{(\da+1)^{\nu+1}}\gs0
\end{equation}
since $a_{0,\ldots,0}=P(0,\ldots,0)>0$.

Define
\begin{equation}\label{3.7}M:=\sum_{j=0}^{(\da+1)^\nu}m_j\mathcal B^j,\end{equation}
where
$$m_j=\begin{cases} \mathcal B-b&\t{if}\ j=(\da+1)^i\ \t{for some}\ i=1,\ldots,\nu,\\\mathcal B-1&\t{otherwise}.\end{cases}$$
Then
$$0\ls M\ls (\mathcal B-1)\sum_{j=0}^{(\da+1)^\nu}\mathcal B^j<N_0:=\mathcal B^{(\da+1)^{\nu}+1}.$$
Let $N_1:=p^2\mathcal B^{(2\da+1)(\da+1)^{\nu}+1}$. Then
$$0\ls (\mathcal B-p)\mathcal B^{(\da+1)^{\nu+1}}<\mathcal B^{(\da+1)^{\nu+1}+1}\ls N_1$$
and hence
$$0\ls T:=M+(\mathcal B-p)\mathcal B^{(\da+1)^{\nu+1}}N_0\ls N_0-1+(N_1-1)N_0<N,$$
where
\begin{equation}\label{3.8}N:=N_0N_1=p^2\mathcal B^{2(\da+1)^{\nu+1}+2}\eq b^{2\da((\da+1)^{\nu+1}+1)}\eq1\pmod{p^2-1}.\end{equation}
(Note that $b\eq1\pmod{p^2-1}$ by (\ref{3.3}).)

Define
\begin{equation}\label{3.9}\mathcal C:=b\mathcal B^{(\da+1)^{\nu}}=b(p^{\al p}b^{\da})^{(\da+1)^{\nu}}.\end{equation}
Let $g\in[0,c\mathcal B^{(\da+1)^\nu})$ and set
$$J:=p(1+a\mathcal B+g)^\da D(\mathcal B)+\sum_{i=0}^{(2\da+1)(\da+1)^\nu}\mathcal B^{i+1}.$$
As $c=(\nu+1)b\ls \mathcal B$ we have $g<N_0$. Note also that $a\mathcal B+(g+1)\ls(a+c)\mathcal B^{(\da+1)^{\nu}}$.
With the aid of (\ref{2.5}) and (\ref{3.5}), we have
\begin{align*}0\ls J\ls&p(a+c)^\da \mathcal B^{\da(\da+1)^\nu}\times\da!L\sum_{i=0}^{(\da+1)^{\nu+1}}\mathcal B^i+\sum_{i=0}^{(2\da+1)(\da+1)^\nu}\mathcal B^{i+1}
\\\ls&(\mathcal B-p)\mathcal B^{\da(\da+1)^\nu}\f{\mathcal B^{(\da+1)^{\nu+1}+1}-1}{\mathcal B-1}+\f{\mathcal B(\mathcal B^{(2\da+1)(\da+1)^{\nu}+1}-1)}{\mathcal B-1}
\\<&\f{(\mathcal B-p)+\mathcal B}{\mathcal B-1}\mathcal B^{(2\da+1)(\da+1)^\nu+1}\ls 2\mathcal B^{(2\da+1)(\da+1)^\nu+1}\ls N_1,
\end{align*}
and hence $0\ls S:=g+JN_0<N_0N_1=N$.
Define
$$R:=(S+T+1)N+T+1,\ X:=\f{N-1}{p-1}R\ \ \t{and}\ \ Y:=N^2.$$
In view of (3.8), we have $p+1\mid X$. Clearly,
$R\gs N+1>\mathcal B=\beta b^\da\gs b>0$, $X\gs\f{p^2-1}{p-1}b\gs3b$, and
$$Y\gs N=p^2(\beta b^\da)^{2(\da+1)^{\nu+1}+2}\gs \max\{b,p^{4p}\}$$
since $\beta\gs p^p$. Thus (3.4) holds.

Below we assume further that $b\in p\uparrow$. Then $\mathcal B,N\in p\uparrow$.
Note that $2\mathcal C\ls c \mathcal B^{(\da+1)^\nu}$ since $2b\ls (\nu+1)b=c$. When $g\in[0,2\mathcal C)$, in view of the last two paragraphs we have
$$\tau_p(S,T)=0\iff \tau_p(g,M)=0\land\tau_p(J,(\mathcal B-p)\mathcal B^{(\da+1)^{\nu+1}})=0$$
since $N_0$ is a power of $p$, and also
$$\tau_p(S,T)=0\iff Y\mid \bi{pX}{X}$$
by Lemma \ref{Lem2.2}.

In view of Lemma \ref{Lem2.4} and (\ref{3.5}), for $z_1,\ldots,z_{\nu}\in[0,b)$ we have
\begin{align*}P(a,z_1,\ldots,z_{\nu})=0\iff& \tau_p\l(J',\l(\f{\mathcal B}p-1\r)\mathcal B^{(\da+1)^{\nu+1}}\r)=0
\\\iff& \tau_p\l(pJ',(\mathcal B-p)\mathcal B^{(\da+1)^{\nu+1}}\r)=0,
\end{align*}
where
$$J':=\(1+a\mathcal B+\sum_{i=1}^\nu z_i\mathcal B^{(\da+1)^i}\)^\da D(\mathcal B)+\f{\mathcal B}p\sum_{i=0}^{(2\da+1)(\da+1)^\nu}\mathcal B^i.$$
If $P(a,z_1,\ldots,z_{\nu})=0$ with $z_1,\ldots,z_{\nu}\in[0,b)$, then $\max\{z_1,\ldots,z_{\nu}\}>0$ since $P(a,0,\ldots,0)>0$,
and hence
\begin{align*}b\ls\mathcal B\ls&\sum_{i=1}^{\nu}z_i\mathcal B^{(\da+1)^i}\ls \sum_{i=1}^{\nu}(b-1)\mathcal B^{(\da+1)^i}
\\\ls&(b-1)\mathcal B^{(\da+1)^\nu}+(\mathcal B-1)\sum_{j=0}^{(\da+1)^\nu-1}\mathcal B^j<b\mathcal B^{(\da+1)^\nu}=\mathcal C.
\end{align*}

Let $\mathcal G\in\{\mathcal C,2\mathcal C\}$. As $b\ls 2b\ls\mathcal B$,
by the above and Lemma \ref{Lem2.3} we have
\begin{align*}&g\in[0,\mathcal G)\land Y\mid \bi{pX}{X}
\\\iff&g\in[0,\mathcal G)\land\tau_p(g,M)=0\land \tau_p(J,(\mathcal B-p)\mathcal B^{(\da+1)^{\nu+1}})=0
\\\iff&\exists z_1\in[0,b)\ldots\exists z_{\nu}\in[0,b)\(g=\sum_{i=1}^\nu z_i\mathcal B^{(\da+1)^i}\land\tau_p(pJ',(\mathcal B-p)\mathcal B^{(\da+1)^{\nu+1}})=0\)
\\\iff&\exists z_1\in[0,b)\ldots z_{\nu}\in[0,b)\(g=\sum_{i=1}^\nu z_i\mathcal B^{(\da+1)^i}\land P(a,z_1,\ldots,z_{\nu})=0\).
\end{align*}

In view of the above, we have completed the proof of Theorem \ref{Th3.1}. \qed
\medskip

\begin{remark}\label{Rem3.1} In the proof of Theorem \ref{Th3.1}, we use $\mathcal B$ and $\mathcal L$ instead of $B$ and $L$ in Lemmas \ref{Lem2.3} and \ref{Lem2.4}.
This is because we will use $B$ and $L$ in later sections for other purposes.
\end{remark}

\section{Working with Lucas sequences}
\setcounter{equation}{0}
\setcounter{theorem}{0}

\begin{lemma}\label{Lem4.1} {\rm (i)} If $0\ls\theta<1$, then $(1-\theta)^n\gs1-n\theta$ for all $n\in\N$.

{\rm (ii)} If $0\ls\theta\ls1/2$, then $1/(1-\theta)\ls1+2\theta$.

{\rm (iii)} For any integers $A\gs2$ and $n\gs0$, we have
\begin{equation}\label{4.1}u_n(A,1)<u_{n+1}(A,1)\ \ \t{and}\ \ (A-1)^n\ls u_{n+1}(A,1)\ls A^n.
\end{equation}
\end{lemma}
\begin{remark}\label{Rem4.1} Lemma \ref{Lem4.1} is easy.
The first part is well known and it can be easily proved by induction.
Part (ii) can be verified directly. Part (iii) can be found in \cite[Lemmas 4 and 8]{S92a}.
\end{remark}

For any integer $A\gs2$, it is known that the solutions of the Pell equation
$$y^2-(A^2-1)x^2=1\ \ (x,y\in\N)$$
are given by $x=u_n(2A,1)$ and $y=v_n(2A,1)/2$ with $n\in\N$.
In this sense, Lemma \ref{Lem4.1}(iii) with $A$ even also appeared in earlier work
 (see, e.g., \cite[Section 2]{MR}).

\begin{lemma}\label{Lem4.2} Let $A,X\in\Z$. Then
\begin{equation}\label{4.2}(A^2-4)X^2+4\in \square\iff X=u_m(A,1)\ \t{for some}\ m\in\Z.
\end{equation}
\end{lemma}
\Proof.  In view of (\ref{1.20}), we have
$$\{u_m(-A,1):\ m\in\Z\}=\{u_m(A,1):\ m\in\Z\}=\{\pm u_n(A,1):\ n\in\N\}.$$
Without any loss of generality, we may simply assume that $A\gs0$.

 If $A\gs2$,  then by \cite[Lemma 9]{S92a} we have
$$X\in\N\land(A^2-4)X^2+4\in\square\iff X=u_n(A,1)\ \t{for some}\ n\in\N,$$
which implies (\ref{4.2}).

For each $n\in\N$, we can easily see that
$$u_n(0,1)=\begin{cases}1&\t{if}\ n\eq1\pmod4,\\0&\t{if}\ n\eq0\pmod2,\\-1&\t{if}\ n\eq-1\pmod4,\end{cases}$$
and
$$ u_n(1,1)=\begin{cases}1&\t{if}\ n\eq1,2\pmod6,\\0&\t{if}\ n\eq0\pmod3,\\-1&\t{if}\ n\eq-1,-2\pmod6.
\end{cases}$$
Therefore, for $A\in\{0,1\}$ we have
$$(A^2-4)X^2+4\in\square\iff X\in\{0,\pm1\}\iff X\in\{u_m(A,1):\ m\in\Z\}.$$

The proof of Lemma \ref{Lem4.2} is now completed. \qed

\begin{lemma} [{\rm See \cite[Theorem 1]{S92a}}]\label{Lem4.3} Let $A,B,C\in\Z$ with $A>1$ and $B\gs0$. Then
\begin{equation}\label{4.3}C=u_B(A,1)\iff C\gs B\land\exists x>0\exists y>0(DFI\in\square),
\end{equation}
where
\begin{equation}\label{4.4}\begin{aligned} &D=(A^2-4)C^2+4,\ E=C^2Dx,\ F=4(A^2-4)E^2+1,
\\&G=1+CDF-2(A+2)(A-2)^2E^2,\ H=C+BF+(2y-1)CF,
\\&I=(G^2-1)H^2+1.
\end{aligned}\end{equation}
Moreover, if $C=u_B(A,1)$ with $B>0$, then for any $Z\in\Z^+$ there are integers $x\gs Z$ and $y\gs Z$ such that $DFI$ is a square.
\end{lemma}
\begin{remark}\label{Rem4.2} Lemma \ref{Lem4.2} (with $A$ not necessarily even) is an extension of Matiyasevich and Robinson's work in \cite[Section 3]{MR}.
The innovation is that we may require arbitrary large solutions when $C=u_B(A,1)$ with $A>1$ and $B\gs1$.
\end{remark}

\begin{lemma} [{\rm See \cite[Theorem 2]{S92a}}]\label{Lem4.4} Let $A,B,C\in\Z$ with $1<|B|<|A|/2-1$. Then
\begin{equation}\label{4.5}C=u_B(A,1)\iff (A-2\mid C-B)\land\exists x\not=0\exists y(DFI\in\square),
\end{equation}
where we adopt the notation in $(\ref{4.4})$.
\end{lemma}
\begin{remark}\label{Rem4.3} This lemma involving integer variables laid the first stone for our proofs of Theorems \ref{Th1.1} and \ref{Th1.2}.
\end{remark}

\begin{lemma}\label{Lem4.5} Let $A,B,U,V\in\Z$ with $B>0$. Then
\begin{equation}\label{4.6}(UV)^{B-1}u_B(A,1)\eq\sum_{r=0}^{B-1}U^{2r}V^{2(B-1-r)}\pmod{U^2-AUV+V^2}.
\end{equation}
\end{lemma}
\Proof. As $u_1(A,1)=1$ and $u_2(A,1)=A$, it is easy to verify (\ref{4.6}) for $B=1,2$.

Below we let $B>2$ and assume that (\ref{4.6}) holds with $B$ replaced by any smaller positive integer. Then
\begin{align*}(UV)^{B-1}u_B(A,1)=&AUV(UV)^{B-2}u_{B-1}(A,1)-U^2V^2(UV)^{B-3}u_{B-2}(A,1)
\\\eq&AUV\sum_{i=0}^{B-2}U^{2i}V^{2(B-2-i)}-U^2V^2\sum_{j=0}^{B-3}U^{2j}V^{2(B-3-j)}
\\\eq&(U^2+V^2)\sum_{i=0}^{B-2}U^{2i}V^{2(B-2-i)}-\sum_{j=0}^{B-3}U^{2j+2}V^{2(B-2-j)}
\\=&\sum_{r=0}^{B-1}U^{2r}V^{2(B-1-r)}
\pmod{U^2-AUV+V^2}.
\end{align*}
This concludes the induction proof of (\ref{4.6}). \qed
\medskip

\begin{remark}\label{Rem4.4} Lemma \ref{Lem4.5} with $U=1$ and $2\mid A$ was first pointed out by Robinson (see also \cite[Lemma 2.22]{J82}) who used it to give a Diophantine representation of the exponential relation
with natural number unknowns.
\end{remark}

\begin{lemma} [{\rm See \cite[Lemma 14]{S92a}}]\label{Lem4.6} Let $B,V$ and $W$ be integers with $B>0$ and $|V|>1$. Then $W=V^B$ if there are $A,C\in\Z$ for which
$|A|\gs\max\{V^{4B},W^4\}$, $C=u_B(A,1)$ and
\begin{equation}\label{4.7} (V^2-1)WC\eq V(W^2-1)\pmod{AV-V^2-1}.\end{equation}
\end{lemma}
\begin{remark}\label{Rem4.5} $A$, $V$ and $W$ Lemma \ref{Lem4.6} are not necessarily positive, they might be negative. In his 1992 PhD thesis \cite{S-PhD}, the author also proved that for $B,V,W\in\Z$ with $B>0$ and $|V|>1$, the equality $W=V^B$ holds if and only if there are integers $A$ and $C$ for which
$|A|\gs\max\{V^{2B},W^2\}$, $C=u_{2B+1}(A,1)$ and
$$ (V-1)WC\eq VW^2-1\pmod{(A^2-2)V-V^2-1}.$$
\end{remark}

The next theorem is motivated by \cite[Lemma 2.25]{J82} on Diophantine representations involving powers of two and central binomial coefficients.
We deal with Diophantine representations involving powers of any prime $p$ and more general binomial coefficients
by only using large variables.

\begin{theorem}\label{Th4.1} Let $p$ be a prime, and let $b\in p\uparrow$ and $g\in\Z^+$. Let $P,Q,X$ and $Y$ be integers with $P>Q>0$ and $X,Y\gs b$.
Suppose that $Y\mid\bi{PX}{QX}$. Then there are integers $h,k,l,w,x,y\gs b$ for which
\begin{equation}\label{4.8}DFI\in\square,\ (U^{2P}V^2-4)K^2+4\in\square,\ pA-p^2-1\mid (p^2-1)WC-p(W^2-1),\end{equation}
\begin{equation}\label{4.9}bw=p^B\ \ \t{and}\ \ 16g^2(C-KL)^2<K^2,\end{equation}
where
\begin{equation}\label{4.10}\begin{aligned} &L:=lY,\ U:=PLX,\ V:=4gwY,\ W:=bw, \ K:=QX+1+k(U^PV-2),
\\& A:=U^Q(V+1),\ B:=PX+1,\ C:=B+(A-2)h,
\end{aligned}\end{equation}
and $D,F,I$ are given by $(\ref{4.4})$.
\end{theorem}
\Proof. Since $b\in p\uparrow$ and
$$p^B\gs p^{PX}\gs(2^X)^P\gs X^2\gs b^2\gs b,$$
we have
$w:=p^B/b\in p\uparrow$ and
\begin{equation}\label{4.11}0<b\ls w\ls W=bw=p^B=p^{PX+1}.\end{equation}
Note that
$$b\ls Y\ls \bi{PX}{QX}\ls\sum_{i=0}^{PX}\bi{PX}i=2^{PX}$$
and
\begin{equation}\label{4.12}8gp^{PX}\ls4gp^B=4gwb\ls V=4gwY\ls 4gWY\ls 4gp^{PX+1}2^{PX}.
\end{equation}
For
\begin{equation}\label{4.13}
\rho:=\f{(V+1)^{PX}}{V^{QX}},
\end{equation}
by the binomial theorem we have
\begin{equation}\label{4.14}\rho
=\f1V\sum_{i=0}^{QX-1}\bi{PX}i\f1{V^{QX-1-i}}+\bi{PX}{QX}+V\sum_{i=QX+1}^{PX}\bi {PX}iV^{i-QX-1}.\end{equation}
As
$$0\ls \f1V\sum_{i=0}^{QX-1}\bi{PX}i\f1{V^{QX-1-i}}<\f1V\sum_{i=0}^{QX}\bi{PX}i\ls\f{2^{PX}}V\ls\f1{8g}<1$$
by (\ref{4.12}), from (\ref{4.14}) we see that
\begin{equation}\label{4.15}\{\rho\}<\f1{8g}\ \ \t{and}\ \lfloor\rho\rfloor=\bi{PX}{QX}+V\sum_{i=QX+1}^{PX}\bi{PX}iV^{i-QX-1}\gs V,\end{equation}
where $\{\rho\}$ is the fractional part of $\rho$, and $\lfloor \rho\rfloor$ is the integral part of $\rho$.
Since $Y$ divides both $\bi{PX}{QX}$ and $V$, we have $l:=\lfloor\rho\rfloor/Y\in\Z$ by (\ref{4.15}). Note that
$$(V+1)^{PX}\gs\rho\gs l=\f{\lfloor\rho\rfloor}Y\gs \f VY=4gw\gs w\gs b$$
and
$$0<U=PLX=\lfloor\rho\rfloor PX\ls\rho PX\ls PX(V+1)^{PX}.$$
Since $A=U^Q(V+1)\gs V+1>2$, by Lemma \ref{Lem4.1}(iii) we have $u_{m+1}(A,1)>u_{m}(A,1)$ for all $m\in\N$.
Clearly, $B=PX+1\gs2X+1\gs3$. Therefore
$$u_B(A,1)\gs u_3(A,1)+(B-3)=A^2-1+B-3=B+(A-2)(A+2).$$
Note that
$$u_B(A,1)\eq u_B(2,1)=B\pmod{A-2}.$$
Thus, for some integer $h\gs A+2$ we have $C=B+(A-2)h=u_B(A,1)$.
Clearly, $A+2\gs V\gs w\gs b$ and hence $h\gs b$. Since $A>1$ and $B>0$, by Lemma \ref{Lem4.3} there are integers $x,y\gs b$
such that $DFI\in\square$.

As
$$u_{QX+1}(U^PV,1)\eq u_{QX+1}(2,1)=QX+1\pmod{U^PV-2},$$
for some $k\in\Z$ we have
\begin{equation}\label{4.16}K=QX+1+k(U^PV-2)=u_{QX+1}(U^PV,1)
\end{equation} and hence
$$(U^{2P}V^2-4)K^2+4\in\square$$
by Lemma \ref{Lem4.2}.
In view of (\ref{4.15}), $U=PLX\gs 2L=2\lfloor\rho\rfloor\gs2V$ and hence
\begin{equation}\label{4.17}U^PV-2\gs 2V-2\gs V\gs w\gs b>0.\end{equation}
If $QX=1$, then $b=1$ since $X\gs b>0$, hence
$$U^PV=u_{2}(U^PV,1)=u_{QX+1}(U^PV,1)=QX+1+k(U^PV-2)=2+k(U^PV-2)$$
and thus $k=1=b$ due to (\ref{4.17}). When $QX>1$, by Lemma \ref{Lem4.1}(iii) and (\ref{4.16})-(\ref{4.17}), we have
\begin{align*}K=u_{QX+1}(U^PV,1)\gs& (U^PV-1)^{QX}=(1+(U^PV-2))^{QX}
\\\gs&1+QX(U^PV-2)+(U^PV-2)^{QX}
\\\gs&1+QX+(U^PV-2)^2\gs 1+QX+b(U^PV-2)
\end{align*}
and hence $k\gs b$.

In light of Lemma \ref{Lem4.5},
\begin{align*}(p^2-1)WC=&p(p^2-1)p^{B-1}u_B(A,1)
\\\eq& p(p^{2B}-1)=p(W^2-1)\pmod{pA-p^2-1}.
\end{align*}

In view of (\ref{4.15})-(\ref{4.17}), $K\gs k\gs b>0$ and
$$A=U^Q(V+1)\gs U^PV\gs U=PLX=\lfloor\rho\rfloor PX\gs VPX\gs 2QX.$$
With the aid of Lemma \ref{Lem4.1}, we have
\begin{align*} \rho\l(1-\f{PX}{U^Q(V+1)}\r)\ls&\rho\l(1-\f1{U^Q(V+1)}\r)^{PX}=\f{(U^Q(V+1)-1)^{PX}}{(U^PV)^{QX}}
\\\ls&\f CK=\f{u_{PX+1}(A,1)}{u_{QX+1}(U^PV,1)}
\\\ls&\f{(U^Q(V+1))^{PX}}{(U^PV-1)^{QX}}=\rho\l(1-\f1{U^PV}\r)^{-QX}
\\\ls&\rho\l(1-\f{QX}{U^PV}\r)^{-1}\ls\rho\l(1+\f{2QX}{U^PV}\r).
\end{align*}
Thus
\begin{equation}\label{4.18}-\rho\f{PX}{U^Q(V+1)}\ls\f CK-\rho\ls \rho\f{2QX}{U^PV},\end{equation}
and hence
$$\l|\f{C}K-\rho\r|\ls\f{2PX}{UV}\rho=\f{2\rho}{LV}=\f{\rho}{\lfloor\rho\rfloor}\times\f2V\ls\f 4V\ls\f1{8g}$$
since $V/(4g)=wY\gs wb=p^B\gs2^{2X+1}\gs8$.
Therefore, in view of (\ref{4.15}), we have
$$\l|\f CK-L\r|=\l|\f CK-\lfloor\rho\rfloor\r|\ls\l|\f CK-\rho\r|+|\rho-\lfloor\rho\rfloor|< \f1{8g}+\f1{8g}=\f1{4g}$$
and hence the inequality in (\ref{4.9}) holds.

Combining the above, we have completed the proof of Theorem \ref{Th4.1}. \qed

The following theorem involving integer variables plays a central role in our later proofs of Theorems \ref{Th1.1} and \ref{Th1.2}.

\begin{theorem}\label{Th4.2} Let $p$ be a prime, and let $b\in\N$ and $g\in\Z^+$. Let $P,Q,X$ and $Y$ be integers with
\begin{equation}\label{4.19}P>Q>0,\ X\gs3b\ \t{and}\ Y\gs\max\{b,p^{4P}\}.\end{equation}
Suppose that there are integers $h,k,l,w,x,y$ with $lx\not=0$ such that both $(\ref{4.8})$ and the inequality
\begin{equation}\label{4.20}4(C-KL)^2<K^2\end{equation} hold, where we adopt the notations in $(\ref{4.10})$ and $(\ref{4.4})$.
Then
\begin{equation}\label{4.21}b\in p\uparrow\ \ \t{and}\ \ Y\mid\bi{PX}{QX}.\end{equation}
\end{theorem}
\Proof. Assume that $W=0$. Then $pA-p^2-1$ divides $(p^2-1)WC-p(W^2-1)=p$ by (\ref{4.8}).
As $p$ is prime and $pA-p^2-1$ is relatively prime to $p$, we must have $pA-p^2-1\in\{\pm1\}$.
Thus $A=p$ or $A=p+1=3$. Note that $U^Q(V+1)=A\gs2$ and $X\gs1$
(since $PLX=U\not=0$ and $X\gs b\gs0$). Hence $|U|=P|L|X\gs2YX\gs2Y\gs2p>A$, which leads to a contradiction since $V+1\eq1\pmod 4$.

By the above, $bw=W\not=0$. Thus $X\gs3b\gs3$ and $PX\gs2\times3b\gs6$. Clearly, $Y\gs p^4\gs4$ by (\ref{4.19}),
\begin{equation}\label{4.22}|A|=|U^Q(V+1)|\gs|U|=PX|L|\gs PXY\gs 4PX>2PX+4\end{equation}
and hence $|A|/2-1>B=PX+1>1$. Recall that $x\not=0$. Also,  $DFI\in\square$ by (\ref{4.8}), and $A-2\mid C-B$ by (\ref{4.10}).
Applying Lemma \ref{Lem4.4} we obtain $C=u_B(A,1)$. In view of (\ref{1.20}) and Lemma \ref{Lem4.1}(iii),
\begin{equation}\label{4.23}|C|=u_B(|A|,1)\ls|A|^{B-1}=|U^Q(V+1)|^{PX}\ls|U^P|^{QX}(|U^PV|-1|)^{PX}
\end{equation}
since $V=4gwY\not=0$ and $|U^PV|-1\gs2|V|-1\gs|V|+1\gs|V+1|$.

As $(U^{2P}V^2-4)K^2+4\in\square$, by Lemma \ref{Lem4.2} we have $K=u_R(U^PV,1)$ for some $R\in\Z$.
Clearly, $(P-Q)X\gs X>2$,
\begin{equation}\label{4.24}
|U^PV|\gs|U|\gs PXY\gs3PX>PX+2QX+4>2.\end{equation}
and
$$QX+1\eq K=u_R(U^PV,1)\eq u_R(2,1)=R\pmod{U^PV-2}.$$

Write $R=QX+1+r(U^PV-2)$ with $r\in\Z$. Suppose that $r\not=0$. By (\ref{4.24}),
$$|R|\gs|r|\times|U^PV-2|-|QX+1|\gs |U^PV|-2-(QX+1)>PX+QX$$
and hence
$$|K|=|u_R(U^PV,1)|=u_{|R|}(|U^PV|,1)\gs(|U^PV|-1)^{|R|-1}\gs(|U^PV|-1)^{PX+QX}$$
with the aid of Lemma \ref{4.1}(iii).
Combining this with (\ref{4.23}) and noting that $|U^PV|\gs4|U^P|>2|U^P|+1$,  we immediately get
$$\l|\f CK\r|\ls \l(\f{|U^P|}{|U^PV|-1}\r)^{QX}<\l(\f12\r)^{QX}\ls\f12.$$
This, together with (\ref{4.20}), yields that
$$|L|\ls\l|L-\f CK\r|+\l|\f CK\r|<\f12+\f12\ls 1,$$
which contradicts $L=lY\not=0$.

By the last paragraph, $R=QX+1$ and hence $K=u_{QX+1}(U^PV,1)$.
As
\begin{equation}\label{4.25}\min\{|A|,|U^PV|\}\gs |U|\gs4PX\gs 4QX
\end{equation} and
$$\l|\f CK\r|=\f{u_{PX+1}(|A|,1)}{u_{QX+1}(|U^PV|,1)},$$
we have
\begin{equation}\label{4.26}
-|\rho|\f{PX}{|U^Q(V+1)|}\ls\l|\f CK\r|-|\rho|\ls |\rho|\f{2QX}{|U^PV|}\end{equation}
in the spirit of the proof of (\ref{4.18}), where $\rho=(V+1)^{PX}/V^{QX}$.
From (\ref{4.25}) and (\ref{4.26}) we deduce that
\begin{equation}\label{4.27}\l|\f CK\r|\gs\f{|\rho|}2.\end{equation}

Note that $|V|\gs 4Y\gs 4p^{4P}\gs4P>4Q$. With the help of Lemma \ref{Lem4.1}(i),
\begin{align*}\f{|V+1|^{Q+1}}{|V|^Q}\gs&\f{(|V|-1)^{Q+1}}{|V|^Q}=(|V|-1)\l(1-\f1{|V|}\r)^Q
\\\gs&(|V|-1)\l(1-\f Q{|V|}\r)\gs(|V|-1)\l(1-\f14\r)\gs\f{|V|-1}2
\end{align*}
and hence
\begin{equation}\label{4.28}
|\rho|\gs\l(\f{|V+1|^{Q+1}}{|V|^Q}\r)^X\gs\l(\f{|V|-1}2\r)^X\gs\l(\f{4Q}2\r)^X\gs2^X\gs 2.
\end{equation}
Combining (\ref{4.20}), (\ref{4.27}) and (\ref{4.28}) we obtain
\begin{equation}\label{4.29}
|L|>\l|\f CK\r|-\f12\gs\f{|\rho|}2-\f12\gs\f{|\rho|}{4}\gs\f1{4}\l(\f{|V|-1}2\r)^X\end{equation}
and hence
\begin{equation}\label{4.30}
|A|\gs|U(V+1)|\gs PX|L|(|V|-1)\gs\f{PX}{2}\l(\f{|V|-1}2\r)^{X+1}\gs\l(\f{|V|-1}2\r)^{X+1}.
\end{equation}
As $|V|-1\gs4Y-1\gs2Y$, from (\ref{4.30}) and (\ref{4.19}) we get
$$|A|\gs Y^{X+1}\gs(p^{4P})^{X+1}\gs p^{4(PX+1)}=p^{4B}.$$
Since
$$\f{|V|-1}2=2|gwY|-\f12\gs|gwY|\gs|wb|=|W|,$$
by (\ref{4.30}) we also have
$|A|\gs |W|^{X+1}\gs W^4$
since $X\gs 3b\gs3$.
As $C=u_B(A,1)$ and
$$(p^2-1)WC\eq p(W^2-1)\pmod{pA-p^2-1}$$
by (\ref{4.8}), applying Lemma \ref{Lem4.6} we obtain $W=p^B$ and thus $bw=p^{PX+1}$.
As $b>0$, we must have $b,w\in p\uparrow$.

Now,
\begin{equation}\label{4.31}
V=4gwY\gs4gwb=4gW\gs 4W=4p^{PX+1}\gs8\times2^{PX}
\end{equation} and hence
$$0\ls \f1V\sum_{i=0}^{QX-1}\f{\bi{PX}i}{V^{QX-1-i}}<\f1V\sum_{i=0}^{PX}\bi{PX}i=\f{2^{PX}}V\ls\f18.$$
Combining this with (\ref{4.14}) we see that
$$\{\rho\}<\f18\ \t{and}\ \lfloor\rho\rfloor=\bi{PX}{QX}+V\sum_{i=QX+1}^{PX}\bi{PX}iV^{i-QX-1}.$$
As $Y$ divides both $L$ and $V$, we have $Y\mid\bi{PX}{QX}$ provided $\lfloor\rho\rfloor=L$.
If
\begin{equation}\label{4.32}\l|\f CK-\rho\r|<\f14,
\end{equation} then
$$|\lfloor\rho\rfloor-L|\ls|\lfloor\rho\rfloor-\rho|+\l|\rho-\f CK\r|+\l|\f CK-L\r|<\f18+\f14+\f12<1$$
with the aid of (\ref{4.20}).
So it suffices to show (\ref{4.32}).

 By (\ref{1.20}),
$$(-A)^{PX}u_{PX+1}(-A,1)=A^{PX}(-1)^{PX}u_{PX+1}(-A,1)=A^{PX}u_{PX+1}(A,1).$$
Thus, in view of Lemma \ref{Lem4.1}(iii), we have
\begin{equation}\label{4.33}
A^{PX}C=A^{PX}u_{PX+1}(A,1)=|A|^{PX}u_{PX+1}(|A|,1)>0\end{equation}
since $|A|=|U^Q(V+1)|\gs V>2$.
Similarly,
\begin{equation}\label{4.34}
(U^PV)^{QX}K=(U^PV)^{QX}u_{QX+1}(U^PV,1)>0.\end{equation}
Now that
$$A^{PX}(U^PV)^{QX}=U^{2PQX}(V+1)^{PX}V^{QX}>0,$$
we must have $CK>0$ by (\ref{4.33}) and (\ref{4.34}). In light of (\ref{4.26}), (\ref{4.29}) and (\ref{4.31}), we finally get
$$\l|\f CK-\rho\r|\ls\rho\f{2PX}{|U|V}=\f{\rho}{|L|}\times\f2V<\f{8}{V}\ls\f1{2^{PX}}\ls\f14.$$
This shows the desired (\ref{4.32}) and thus concludes our proof of Theorem \ref{Th4.2}. \qed

\section{Proofs of Theorem \ref{Th1.1} and Corollary \ref{Cor1.2}}
\setcounter{equation}{0}
\setcounter{theorem}{0}

During their reduction of unknowns in Diophantine representations, Matiyasevich and Robinson \cite{MR}
introduced for each $k\in\Z^+$ the polynomial
\begin{equation}\label{5.1}J_k(x_1,\ldots,x_k,x):=\prod_{\ve_1,\ldots,\ve_k\in\{\pm1\}}
\l(x+\ve_1\sqrt{x_1}+\ve_2\sqrt{x_2}X+\ldots+\ve_k\sqrt{x_k}X^{k-1}\r)\end{equation}
with $X=1+\sum_{i=1}^kx_i^2$. They showed that this polynomial has integer coefficients and that $A_1,\ldots,A_k\in\Z$
are all squares if and only if $J_k(A_1,\ldots,A_k,x)=0$ for some $x\in\Z$.

\begin{lemma}[See \cite{MR}]\label{Lem5.1}  Let $A_1,\ldots,A_k,R,S$ and $T$ be integers with $S\not=0$. Then
\begin{equation}\label{5.2}\begin{aligned} &A_1\in\square\land\ldots\land A_k\in\square\land S\mid T\land R>0
\\\iff&\exists n\gs0[M_k(A_1,\ldots,A_k,S,T,R,n)=0],
\end{aligned}\end{equation}
where
\begin{align*}&M_k(x_1,\ldots,x_k,w,x,y,z)
\\=&\prod_{\ve_1,\ldots,\ve_k\in\{\pm1\}}\(x^2+w^2z-w^2(2y-1)\(x^2+X^k+\sum_{j=1}^k\ve_j\sqrt{x_j}X^{j-1}\)\)
\\=&(w^2(1-2y))^{2^k}J_k\bigg(x_1,\ldots,x_k,x^2+X^k+\f{x^2+w^2z}{w^2(1-2y)}\bigg)\in\Z[x_1,\ldots,x_k,w,x,y,z]
\end{align*}
with $X=1+\sum_{j=1}^k x_j^2$.
\end{lemma}
\begin{remark}\label{Rem5.1} If $A_1,\ldots,A_k\in\square$, and $R,S,T$ are integers with $R>0$, $S\not=0$ and $S\mid T$, then we can easily see that
$$M_k(A_1,\ldots,A_k,S,T,R,m)=0,$$
where
$$m=(2R-1)(T^2+X^k+\sqrt{A_1}X^0+\ldots+\sqrt{A_k}X^{k-1})-\f {T^2}{S^2}\gs X\gs\max\{A_1,\ldots,A_k\}$$
with $X=1+\sum_{j=1}^kA_j^2$.
\end{remark}

\begin{lemma}\label{Lem5.2} For any $A_1,\ldots,A_k,S,T\in\Z$ with $S\not=0$, we have
\begin{equation}\label{5.3}\begin{aligned} &A_1\in\square\land\cdots\land A_k\in\square\land S\mid T
\iff\exists z[H_k(A_1,\ldots,A_k,S,T,z)=0],
\end{aligned}\end{equation}
where
\begin{equation}\label{5.4}
H_k(x_1,\ldots,x_k,x,y,z):=x^{2^k}J_k\l(x_1,\ldots,x_k,z-\f yx\r)\in\Z[x_1,\ldots,x_k,x,y,z].
\end{equation}
\end{lemma}
\begin{remark}\label{Rem5.2} This is \cite[Lemma 17]{S92a} motivated by Lemma \ref{Lem5.1}. Note that $z$ in (\ref{5.3}) is an integer variable.
\end{remark}

\begin{lemma}\label{Lem5.3} Let $m\in\Z$. Then
\begin{equation}\label{5.5}m\gs0\iff\exists x\not=0[(3m-1)x^2+1\in\square].\end{equation}
\end{lemma}
\Proof. Clearly, $(3\times0-1)1^2+1\in\square$. If $m<0$ and $x\in\Z\sm\{0\}$, then $(3m-1)x^2+1\ls-4+1<0$.
If $m>0$, then $3m-1>0$ and $3m-1\not\in\square$, hence the Pell equation $y^2-(3m-1)x^2=1$
has infinitely many integral solutions and thus $(3m-1)x^2+1\in\square$ for some nonzero integer $x$.
Thus (\ref{5.5}) always holds.
\qed

\medskip
\noindent{\it Proof of Theorem \ref{Th1.1}}. By Matiyasevich's theorem, $\mathcal A$ is a Diophantine set. Let $p$ be a prime. Then (\ref{3.1}) and (\ref{3.2}) hold with $b,\mathcal C$ and $X,Y\in\Z[a,f,g]$ as in Theorem \ref{Th3.1}. Set $P=p$ and $Q=1$, and adopt the notations in (\ref{4.4}) and (\ref{4.10}).

(i) Suppose that $a\in\mathcal A$. By (\ref{3.1}), for any $Z\in\Z^+$ we may take $f\gs Z$ with $b\in\square$ and $b\in p\uparrow$,
and $g\in[b,\mathcal C)$ with $Y$ dividing $\bi{PX}{QX}=\bi{pX}X$. Clearly,
$$0<f\ls b\ls g<\mathcal C<2\mathcal C.$$
As (\ref{3.4}) is valid, by Theorem \ref{Th4.1} there are integers $h,k,l,w,x,y\gs b$
such that both (\ref{4.8}) and (\ref{4.9}) hold.
Thus
$$4(C-KL)^2+\f{g^2K^2}{8\mathcal C^3}<\f{K^2}{4g^2}+\f{K^2}{8g}\ls\f{K^2}g$$
and hence
\begin{equation}\label{5.6}
O:=f^2l^2x^2(8\mathcal C^3gK^2-g^2(32(C-KL)^2\mathcal C^3+g^2K^2))>0.\end{equation}
Note that $g,h,k,l,w,x,y\gs b\gs f\gs Z.$
In view of (\ref{4.8}) and the facts $b\in\square$ and $O>0$, by Remark \ref{Rem5.1} we have
\begin{equation}\label{5.7}
P_{\mathcal A}(a,f,g,h,k,l,w,x,y,m)=0\end{equation}
for some integer $m\gs b\gs f\gs Z$, where
\begin{equation}\label{5.8}\begin{aligned} &P_{\mathcal A}(a,f,g,h,k,l,w,x,y,m)
\\=&M_3(b,DFI,(U^{2P}V^2-4)K^2+4,
\\&pA-p^2-1,(p^2-1)WC-p(W^2-1),O,m).
\end{aligned}\end{equation}
Note that $P_{\mathcal A}(z_0,z_1,\ldots,z_9)\in\Z[z_0,z_1,\ldots,z_9]$. So (\ref{1.3}) has been proved.

Let $a\in\N$, and assume that there are integers $m\gs0$ and $f,g,h,k,l,w,x,y$ satisfying (\ref{5.7}). By Lemma \ref{Lem5.1} we have
(\ref{4.8}), also $b\in\square$ and $O>0$.  By (\ref{5.6}), $fglx\not=0$.
As $b\gs0$ and $f\not=0$, we have $b>0$ and hence $\mathcal C>0$. It follows from (\ref{5.6}) that
$$\f{K^2}g>4(C-KL)^2+\f{g^2K^2}{8\mathcal C^3}\gs\f{g^2K^2}{8\mathcal C^3}\gs0.$$
Thus $K\not=0$ and $0<g<2\mathcal C$. Now, (\ref{3.4}), (\ref{4.19}) and (\ref{4.20}) all hold. By Theorem \ref{Th4.2},
we have $b\in p\uparrow$ and $\bi{pX}{X}=\bi{PX}{QX}\eq0\pmod Y$. Hence $a\in\mathcal A$ by (\ref{3.2}).
This proves (\ref{1.2}).

In view of the above, we have proved the first part of Theorem \ref{Th1.1}.

(ii) By the above,  a nonnegative integer $a$ belongs to $\mathcal A$, if and only if there are integers
$f,g,h,k,l,w,x,y$ such that $b\in\square$, $O>0$, and (\ref{4.8}) holds. By Lemma \ref{Lem5.3},
$$O>0\iff O-1\gs0\iff\exists z\not=0[(3O-4)z^2+1\in\square].$$
In light of Lemma \ref{Lem5.2}, we have
\begin{align*}&b\in\square,\ (3O-4)z^2+1\in\square,\ \t{and (\ref{4.8}) holds}
\\\iff &\exists m[Q_{\mathcal A}(a,f,g,h,k,l,m,w,x,y,z)=0],
\end{align*}
where
\begin{align*}&Q_{\mathcal A}(a,f,g,h,k,l,m,w,x,y,z)
\\=&H_4(b,(3O-4)z^2+1,DFI,(U^{2P}V^2-4)K^2+4,
\\&\quad\ pA-p^2-1,(p^2-1)WC-p(W^2-1),m).
\end{align*}
Note that $Q_{\mathcal A}(z_0,z_1,\ldots,z_{10})\in\Z[z_0,z_1,\ldots,z_{10}]$ and (\ref{1.4}) holds.

The proof of Theorem \ref{Th1.1} is now complete. \qed

\medskip
\noindent{\it Proof of Corollary \ref{Cor1.2}}. Let $\mathcal A\se\N$ be a nonrecursive r.e. set. By Theorem \ref{Th1.1}(i), there is a polynomial
$P_{\mathcal A}(z_0,z_1,\ldots,z_9)\in\Z[z_0,z_1,\ldots,z_9]$ such that for any $a\in\N$ we have
$$ a\in\mathcal A\iff \exists z_1\ldots\exists z_8\exists z_9\gs0[P_{\mathcal A}(a,z_1,\ldots,z_9)=0].$$
Thus, with the aid of Lemma \ref{Lem5.3},
\begin{align*}a\not\in \mathcal A\iff&\neg  \exists z_1\ldots\exists z_8\exists z_9[z_9\gs0\land P_{\mathcal A}(a,z_1,\ldots,z_9)=0]
\\\iff&\forall z_1\ldots\forall z_8\forall z_9[z_9<0\lor P_{\mathcal A}(a,z_1,\ldots,z_9)\not=0]
\\\iff&\forall z_1\ldots\forall z_8\forall z_9[-z_9-1\gs0\lor P_{\mathcal A}(a,z_1,\ldots,z_9)\not=0]
\\\iff&\forall z_1\ldots\forall z_8\forall z_9[\exists x\not=0((3(-z_9-1)-1)x^2+1\in\square)
\\&\ \ \ \lor \exists x\not=0(P_{\mathcal A}(a,z_1,\ldots,z_9)=x)]
\\\iff&\forall z_1\ldots\forall z_8\forall z_9\exists x\not=0[1-(3z_9+4)x^2\in\square
\lor P_{\mathcal A}(a,z_1,\ldots,z_9)=x].
\end{align*}
In view of (\ref{1.5}),
\begin{align*}&\exists x\not=0[1-(3z_9+4)x^2\in\square\lor P_{\mathcal A}(a,z_1,\ldots,z_9)=x]
\\\iff&\exists x\not=0\exists y[(1-(3z_9+4)x^2-y^2)(P_{\mathcal A}(a,z_1,\ldots,z_9)-x)=0]
\\\iff&\exists x_1\exists x_2\exists y[(1-(3z_9+4)(2x_1+1)^2(3x_2+1)^2-y^2)
\\&\ \ \ \times(P_{\mathcal A}(a,z_1,\ldots,z_9)-(2x_1+1)(3x_2+1))=0].
\end{align*}
Therefore $\forall^9\exists^3$ over $\Z$ is undecidable.

By Theorem \ref{Th1.1}(ii), there is a polynomial $Q_{\mathcal A}(z_0,\ldots,z_{10})\in\Z[z_0,\ldots,z_{10}]$
such that (\ref{1.4}) holds for any $a\in\N$.
Hence
\begin{align*} a\not\in\mathcal A\iff&\neg \exists z_1\ldots\exists z_9\exists z_{10}[z_{10}\not=0\land Q_{\mathcal A}(a,z_1,\ldots,z_9,z_{10})=0].
\\\iff&\forall z_1\ldots\forall z_9\forall z_{10}[z_{10}=0\lor Q_{\mathcal A}(a,z_1,\ldots,z_9,z_{10})\not=0]
\\\iff&\forall z_1\ldots\forall z_9\forall z_{10}\exists x\exists y[z_{10}(Q_{\mathcal A}(a,z_1,\ldots,z_9,z_{10})-(2x+1)(3y+1))=0]
\end{align*}
by using (\ref{1.5}). Thus $\forall^{10}\exists^2$ over $\Z$ is undecidable.

So far we have completed the proof of Corollary \ref{Cor1.2}. \qed

\section{Proofs of Theorems \ref{Th1.2} and \ref{Th1.3}}
\setcounter{equation}{0}
\setcounter{theorem}{0}

\begin{lemma}\label{Lem6.1}
{\rm (i)} Any integer can be written as $2^{\da}(x^2-y^2)$ with $\da\in\{0,1\}$ and $x,y\in\Z$.
Also, each integer can be written as $2^{\da}(p_8(x)-p_8(y))$ with $\da\in\{0,1\}$ and $x,y\in\Z$.

{\rm (ii)} Any positive odd integer can be written as $x^2+y^2+2z^2$ with $x,y,z\in\Z$. Also,
each positive odd integer can be written as $p_8(x)+p_8(y)+2p_8(z)$ with $x,y,z\in\Z$.

{\rm (iii)} For any $x\in\Z$, we have $x=T_{x}-T_{-x}=p_5(-x)-p_5(x)$. Also,
\begin{equation}\label{6.1}\{T_x+T_y+T_z:\ x,y,z\in\Z\}=\N=\{p_5(x)+p_5(y)+p_5(z):\ x,y,z\in\Z\}.
\end{equation}
\end{lemma}
\Proof. (i) Clearly, $0=0^2-0^2$.
Write $n\in\Z\sm\{0\}$ as $2^km$ with $k\in\N$, $m\in\Z$ and $2\nmid m$. If $k$ is even, then
$$n=\l(2^{k/2}\f{m+1}2\r)^2- \l(2^{k/2}\f{m-1}2\r)^2.$$
If $k$ is odd, then
$$n=2\l(2^{(k-1)/2}\f{m+1}2\r)^2- 2\l(2^{(k-1)/2}\f{m-1}2\r)^2.$$

Let $n\in\Z$. If $n=4x$ for some $x\in\Z$, then
$n=p_8(-x)-p_8(x)$. If $n=2x+1$ for some $x\in\Z$, then $n=p_8(x+1)-p_8(-x)$.
If $n=2x$ with $x$ odd, then
$$n=2\l(p_8\l(\f{x+1}2\r)-p_8\l(\f{1-x}2\r)\r).$$

In view of the above, we have proved part (i) of Lemma \ref{Lem6.1}.

(ii) The first assertion in part (ii) is well known. Actually, it can be deduced from the Gauss-Legendre theorem on sums of three squares.
For any $n\in\N$, we can write $4n+2$ as $x^2+y^2+(2z)^2$ with $x,y,z\in\Z$ and $x\eq y\pmod 2$, and hence
$$2n+1=\f{x^2+y^2}2+2z^2=\l(\f{x+y}2\r)^2+\l(\f{x-y}2\r)^2+2z^2.$$

Now we prove the second assertion in part (ii). Let $n\in\Z^+$. By \cite[Lemma 4.3(ii)]{S16}, $6n+1=x^2+y^2+2z^2$ for some $x,y,z\in\Z$ with $3\nmid xyz$.
As $x$ or $-x$ is congruent to $-1$ modulo $3$, without loss of generality we may assume that $x=3u-1$ for some $u\in\Z$. Similarly, we may assume
that $y=3v-1$ and $z=3w-1$ for some $v,w\in\Z$. Thus
\begin{align*}6n+1=&(3u-1)^2+(3v-1)^2+2(3w-1)^2
\\=&(3p_8(u)+1)+(3p_8(v)+1)+2(3p_8(w)+1)
\end{align*}
and hence $2n-1=p_8(u)+p_8(v)+2p_8(w)$.

(iii) The first assertion in Lemma \ref{Lem6.1}(iii) can be easily seen.
The first equality in (\ref{6.1}) was conjectured by Fermat and proved by Gauss (see, e.g.,
\cite[p.\,27]{N96}).
The second equality in (\ref{6.1}) was first observed by Guy \cite{G94} (see also the paragraph in \cite{S15} containing \cite[(1.4)]{S15} for a supplement to Guy's proof).
\qed

\medskip
\noindent{\it Proof of Theorem \ref{Th1.2}}.
As
$$8T_z+1=(2z+1)^2,\ 3p_8(z)+1=(3z-1)^2\ \t{and}\ 24p_5(z)+1=(6z-1)^2,$$
we get
\begin{align}\label{6.2}\{8t+1:\ t\in\t{Tri}\}&=\{z^2:\ z\in\Z\land 2\nmid z\},
\\\label{6.3}\{3q+1:\ q\in\t{Octa}\}&=\{z^2:\ z\in\Z\land 3\nmid z\},
\\\label{6.4}\{24r+1:\ r\in\t{Pen}\}&=\{z^2:\ z\in\Z\land 2\nmid z\land 3\nmid z\}.
\end{align}

Let $p$ be a prime. Set $P=p$ and $Q=1$.
$P_{\mathcal A}(a,f,g,h,k,l,w,x,y,m)$ given by (\ref{5.8}) can be written as
$$Q_p(a,f,g,h,k,l,w,x^2,y,m)$$ with $Q_p(z_0,\ldots,z_9)\in\Z[z_0,\ldots,z_9]$.
(Actually, $F,G,H$ and $I$ in (\ref{4.4}) involve $E^2=C^4D^2x^2$.)
When $b\in\square$, $w\in\Z$ and $bw=p^{pX+1}$ with $p+1\mid X$, we have
$$b\f wp=p^{pX}=\l(p^{pX/2}\r)^2\in\square$$
and hence $w=ps$ for some $s\in \square\cap p\uparrow$.
In view of (\ref{3.4}) in Theorem \ref{Th3.1} and (\ref{4.9}) in Theorem \ref{Th4.1}, by modifying the proof of Theorem \ref{Th1.1}(i) slightly
we see that
\begin{equation}\label{6.5}\begin{aligned} a\in\mathcal A\iff &Q_p(a,f,g,h,k,l,ps,x^2,y,m)=0
\\&\ \t{for some}\ f,g,h,k,l,m,s,x,y\in\Z\ \t{with}\ m\gs0\ \t{and}\ s\in\square.
\end{aligned}\end{equation}
Similarly, in view of (\ref{6.2}), when $p\not=2$ we have
\begin{equation}\label{6.6}\begin{aligned} a\in\mathcal A\iff &Q_p(a,f,g,h,k,l,p(8t+1),x^2,y,m)=0
\\&\ \t{for some}\ f,g,h,k,l,m,t,x,y\in\Z\ \t{with}\ m\gs0\ \t{and}\ t\in\t{Tri}.
\end{aligned}\end{equation}
With the help of (\ref{6.3}), if $p\not=3$ then
\begin{equation}\label{6.7}\begin{aligned} a\in\mathcal A\iff &Q_p(a,f,g,h,k,l,p(3q+1),x^2,y,m)=0
\\&\ \t{for some}\ f,g,h,k,l,m,q,x,y\in\Z\ \t{with}\ m\gs0\ \t{and}\ q\in\t{Octa}.
\end{aligned}\end{equation}
In view of (\ref{6.4}), when $p>3$ we have
\begin{equation}\label{6.8}\begin{aligned} a\in\mathcal A\iff &Q_p(a,f,g,h,k,l,p(24r+1),x^2,y,m)=0
\\&\ \t{for some}\ f,g,h,k,l,m,r,x,y\in\Z\ \t{with}\ m\gs0\ \t{and}\ r\in\t{Pen}.
\end{aligned}\end{equation}

When $a\in\mathcal A$, by the proof of Theorem \ref{Th1.1}(i) and the above arguments, there are $f,g,h,k,l,s,x,y\in\Z$ with $s\in\square\cap p\uparrow$
for which
$$b\in p\uparrow,\ b>1,\ b\in\square,\ DFI\in\square,\ (U^{2p}V^2-4)K^2+4\in\square,$$
$$pA-p^2-1\mid(p^2-1)WC-p(W^2-1)\ \t{and}\ O>0$$
with $w=ps$. Since $U=pXL$ is divisible by $p(p+1)$, we see that $U,A,D,(U^{2p}V^2-4)K^2+4$
are all even. If we take $p=2$, then $2\mid b$, $2\nmid pA-p^2-1$, and
$$X_0=1+b^2+(DFI)^2+((U^{2p}V^2-4)K^2+4)^2\eq1\pmod 2;$$
 hence by Remark \ref{Rem5.1} we have
 $$Q_p(a,f,g,h,k,l,ps,x^2,y,m)=0$$
 for some $m\in\N$ with
 \begin{align*}m\eq&(2O-1)(((p^2-1)WC-p(W^2-1))^2+1)
\\&-((p^2-1)WC-p(W^2-1))^2
 \\\eq&1\ (\mo\ 2).
 \end{align*}

In view of parts (i)-(ii) of Lemma \ref{Lem6.1} and the above, by taking $p=2$ we get
 \begin{align*}a\in\mathcal A\iff&Q_p(a,f,g,h,k,l,ps,x,y,m+u+2v)=0
 \\&\t{for some}\ f,g,h,k,l,s,x,y,m,u,v\in\Z\ \t{with}\ s,x,m,u,v\in\square
 \\\iff&\prod_{\da_1,\ldots,\da_6\in\{0,1\}}Q_p(a,2^{\da_1}(f_1-f_2),2^{\da_2}(g_1-g_2),2^{\da_3}(h_1-h_2),
 \\&2^{\da_4}(k_1-k_2),2^{\da_5}(l_1-l_2),ps,x,2^{\da_6}(y_1-y_2),m+u+2v)=0
 \\&\t{for some}\ f_1,f_2,g_1,g_2,h_1,h_2,k_1,k_2,l_1,l_2,s,x,y_1,y_2,m,u,v\in\square
 \end{align*}
and
 \begin{align*}a\in\mathcal A\iff&Q_p(a,f,g,h,k,l,p(3q+1),x^2,y,m+u+2v)=0
 \\&\t{for some}\ f,g,h,k,l,x,y,m,u,v\in\Z\ \t{with}\ q,m,u,v\in\t{Octa}
 \\\iff&\prod_{\da_1,\ldots,\da_7\in\{0,1\}}Q_p(a,2^{\da_1}(f_1-f_2),2^{\da_2}(g_1-g_2),2^{\da_3}(h_1-h_2),2^{\da_4}(k_1-k_2),
 \\&2^{\da_5}(l_1-l_2),p(3q+1),2^{2\da_6}(x_1-x_2)^2,2^{\da_7}(y_1-y_2),m+u+2v)=0
 \\&\t{for some}\ f_1,f_2,g_1,g_2,h_1,h_2,k_1,k_2,l_1,l_2,x_1,x_2,y_1,y_2,q,m,u,v\in\t{Octa}.
 \end{align*}
Similarly, by taking $p>3$ and noting Lemma \ref{Lem6.1}(iii), (\ref{6.6}) and (\ref{6.8}) we obtain
\begin{align*}a\in\mathcal A\iff&Q_p(a,f,g,h,k,l,p(8t+1),x^2,y,m+u+v)=0
 \\&\t{for some}\ f,g,h,k,l,x,y,m,u,v\in\Z\ \t{with}\ t,m,u,v\in\t{Tri}
 \\\iff&Q_p(a,f_1-f_2,g_1-g_2,h_1-h_2,k_1-k_2,
 \\&l_1-l_2,p(8t+1),(x_1-x_2)^2,y_1-y_2,m+u+v)=0
 \\&\t{for some}\ f_1,f_2,g_1,g_2,h_1,h_2,k_1,k_2,l_1,l_2,x_1,x_2,y_1,y_2,t,m,u,v\in\t{Tri}.
 \end{align*}
 and
\begin{align*}a\in\mathcal A\iff&Q_p(a,f,g,h,k,l,p(24r+1),x^2,y,m+u+v)=0
 \\&\t{for some}\ f,g,h,k,l,x,y,m,u,v\in\Z\ \t{with}\ r,m,u,v\in\t{Pen}
 \\\iff&Q_p(a,f_1-f_2,g_1-g_2,h_1-h_2,k_1-k_2,
 \\&l_1-l_2,p(24r+1),(x_1-x_2)^2,y_1-y_2,m+u+v)=0
 \\&\t{for some}\ f_1,f_2,g_1,g_2,h_1,h_2,k_1,k_2,l_1,l_2,x_1,x_2,y_1,y_2,r,m,u,v\in\t{Pen}.
 \end{align*}

The proof of Theorem \ref{Th1.2} is now completed. \qed

\begin{lemma}\label{Lem6.2} {\rm (i) (See Putnam \cite{P69})} For any polynomial $P(x)\in\Z[x]$, we have
\begin{equation}\label{6.9}\N\cap\{(x+1)(1-P(x)^2)-1:\ x\in\N\}=\{x\in\N:\ P(x)=0\}.
\end{equation}

{\rm (ii) (See Sun \cite{S16})} Each $n\in\N$ can be written as the sum of four generalized octagonal numbers, i.e.,
$n=p_8(z_1)+p_8(z_2)+p_8(z_3)+p_8(z_4)$ for some $z_1,z_2,z_3,z_4\in\Z$.
\end{lemma}
\begin{remark}\label{Rem6.1} (\ref{6.9}) is a simple fact which can be easily seen, nevertheless it's a useful trick due to Putnam \cite{P69}.
The author's result Lemma \ref{Lem6.2}(ii) is quite similar to Lagrange's four-square theorem.
\end{remark}

\medskip
\noindent{\it Proof of Theorem \ref{Th1.3}}. (i) For any polynomial $P(z_0,z_1,\ldots,z_n)\in\Z[z_0,z_1,\ldots,z_n]$, we define
$$P^*(z_0,z_1,\ldots,z_n):=(z_0+1)(1-P(z_0,z_1,\ldots,z_n)^2)-1.$$

Let $P_{\mathcal A}$ be as in Theorem \ref{Th1.1}(i). In view of (\ref{1.2}) and (\ref{1.3}),
\begin{equation}\label{6.10}\{a\in\N:\ \exists z_1\ldots\exists z_8\exists z_9\gs0[P_{\mathcal A}(a,z_1,\ldots,z_9)=0]\}=\mathcal A.\end{equation}
Combining (\ref{6.10}), Lemma \ref{Lem6.2}(i) and (\ref{1.1}), we see that
\begin{align*}\mathcal A=&\N\cap\{P_{\mathcal A}^*(z_0,z_1,\ldots,z_8,z_9):\ z_0,z_9\in\N\ \t{and}\ \ z_1,\ldots,z_8\in\Z\}
\\=&\N\cap\{P_{\mathcal A}^*(z_{12}^2+z_{13}^2+z_{14}^2+z_{14},z_1,\ldots,z_8,z_9^2+z_{10}^2+z_{11}^2+z_{11}:\, z_1,\ldots,z_{14}\in\Z\}.
\end{align*}

Let $P_4,P_3,P_5,P_8$ be as in Theorem \ref{Th1.2}. Then
$$\{a\in\N:\ \exists z_1\in\square\ldots\exists z_{17}\in\square[P_4(a,z_1,\ldots,z_{17})=0]\}=\mathcal A.$$
Combining this with Lemma \ref{Lem6.2}(i) and Lagrange's four-square theorem, we obtain
\begin{align*}\mathcal A=&\N\cap\{P_4^*(z_0,z_1,\ldots,z_{17}):\ z_0\in\N\ \t{and} \ z_1,\ldots,z_{17}\in\square\}
\\=&\N\cap\{P_4^*(z_{18}+z_{19}+z_{20}+z_{21},z_1,\ldots,z_{17}):\, z_1,\ldots,z_{21}\in\square\}.
\end{align*}
Similarly, by (\ref{1.12}) and Lemma \ref{Lem6.2}(ii), we have
\begin{align*}\mathcal A=&\N\cap\{P_8^*(z_0,z_1,\ldots,z_{18}):\ z_0\in\N\ \t{and} \ z_1,\ldots,z_{18}\in\t{Octa}\}
\\=&\N\cap\{P_8^*(z_{19}+z_{20}+z_{21}+z_{22},z_1,\ldots,z_{18}):\, z_1,\ldots,z_{22}\in\t{Octa}\}.
\end{align*}
In view of (\ref{1.12}), (\ref{6.1}) and Lemma \ref{Lem6.2}(i), we also have
$$\mathcal A=\N\cap\{P_3^*(z_{19}+z_{20}+z_{21},z_1,\ldots,z_{18}):\, z_1,\ldots,z_{21}\in\t{Tri}\}$$
and
$$\mathcal A=\N\cap\{P_5^*(z_{19}+z_{20}+z_{21},z_1,\ldots,z_{18}):\, z_1,\ldots,z_{21}\in\t{Pen}\}.$$
This concludes the proof of Theorem \ref{Th1.3}(i).

(ii) Clearly, $\mathcal P$ is an r.e. set. Applying Theorem \ref{Th1.2} with $\mathcal A=\mathcal P$, we see that
$$\mathcal P=\{a\in\N:\ \exists z_1\ldots\exists z_{17}[Q(a,z_1^2,\ldots,z_{17}^2)=0]\}$$
for some polynomial $Q(z_0,z_1,\ldots,z_{17})\in\Z[z_0,z_1,\ldots,z_{17}]$. By Lemma \ref{Lem6.1}(ii), any prime can be written as
$x^2+y^2+2z^2$ with $x,y,z\in\Z$. Thus
$$\mathcal P=\{x^2+y^2+2z^2:\ x,y,z\in\Z\land\exists z_1\ldots\exists z_{17}[Q(x^2+y^2+2z^2,z_1^2,\ldots,z_{17}^2)=0]\},$$
and hence (\ref{1.15}) holds with
$$\hat P(z_1,\ldots,z_{20}):=(z_{18}+z_{19}+2z_{20}+1)(1-Q(z_{18}+z_{19}+2z_{20},z_1,\ldots,z_{17})^2)-1.$$
Similarly, by Theorem \ref{Th1.2} and the second assertion in Lemma \ref{Lem6.1}(ii), (\ref{1.16}) holds for certain polynomial $\tilde P(z_1,\ldots,z_{21})
\in\Z[z_1,\ldots,z_{21}]$. (Note that $2=p_8(0)+p_8(0)+2p_8(1)$ and $p_8(-z)=z(3z+2)$ for $z\in\Z$.)

The proof of Theorem \ref{Th1.3} is now completed. \qed
\medskip

\Ack. This work was supported by National Natural Science Foundation of China (Grant No. 11971222).



\begin{thebibliography}{99}


\bibitem {B68} Baker A.  Contributions to the theory of diophantine equations I: On the representation of integers by binary forms. Philos. Trans. Roy. Soc. London (Ser. A), 1968, 263: 173--191

\bibitem {C80} Cutland N.  Computability. Cambridge: Cambridge Univ. Press, 1980

\bibitem {D73} Davis M.  Hilbert's tenth problem is unsolvable. Amer. Math. Monthly, 1973, 80: 233--269

\bibitem {DMR} Davis M, Matiyasevich Yu, Robinson J.  Hilbert's tenth problem.
Diophantine equations: positive aspects of a negative solution. In: Mathematical Decelopments Arising from Hilbert Problems (Proc. Sympos. Pure Math., Vol. 28), Providence, R.I.:  Amer. Math. Soc., 1976, 323--378

\bibitem {DPR} Davis M, Putnam H, Robinson J.  The decision problem for exponential diophantine equations. Ann. of Math., 1961, 74(2): 425--436

\bibitem {De75} Denef J.  Hilbert's Tenth Problem for quadratic rings. Proc. Amer. Math. Soc.,
 1975, 48: 214--220

\bibitem {De78} Denef J.  The Diophantine problem for polynomial rings and fields of rational functions. Trans. Amer. Math. Soc., 1978, 242: 391--399

\bibitem {DL} Denef J, Lipshitz L.  Diophantine sets over some rings of algebraic integers.
J. London Math. Soc., 1978, 18: 385--391

\bibitem {FW} Flath D,  Wagon S.  How to pick out the integers in the rationals: an application of number theory to logic? Amer. Math. Monthly, 1991, 98: 812--823

\bibitem {G94} Guy, R K. Every number is expressible as the sum of how many polygonal numbers?
Amer. Math. Monthly, 1994, 101: 169--172

\bibitem {J81} Jone J P.  Classification of quantifier prefixes over Diophantine equations. Z. Math. Logik Grundlag. Math., 1981, 27: 403--410

\bibitem {J82} Jones J P.  Universal Diophantine equation. J. Symbolic Logic, 1982, 47: 549--571

\bibitem {JM} Jone J P, Matiyasevich Yu.  Register machine proof of the theorem on exponential Diophantine representation of enumerable sets. J. Symbolic Logic, 1984, 49: 818--829

\bibitem {K} Koenigsmann J. Defining $\Z$ in $\Q$. Ann. of Math., 2016, 183: 73--93

\bibitem {M70} Matiyasevich Yu. Enumerable sets are diophantine. Dokl. Akad. Nauk SSSR,
 1970, 191: 279--282;
English translation with addendum, Soviet Math. Doklady, 1970, 11: 354--357

\bibitem {M77}  Matiyasevich Yu.  Some purely mathematical results inspired by mathematical logic.
In: Logic, Foundations of Mathematics and Computability Theory
(London, Ont., 1975). Reidel, Dordrecht, 1977, Part I, 121--127

\bibitem {M81} Matiyasevich Yu.  Primes are nonnegative values of a polynomial in 10 variables.
 J. Soviet Math., 1981, 15: 33--44

\bibitem {M93} Matiyasevich Yu. Hilbert's Tenth Problem. Cambridge, Massachusetts: MIT Press, 1993

\bibitem {MR} Matiyasevich Yu, Robinson J.  Reduction of an arbitrary diophantine equation to one in 13 unknowns. Acta Arith., 1975, 27, 521--553

\bibitem {N96} Nathanson M B.  Additive Number Theory: The
Classical Bases. Grad. Texts in Math., vol. 164,
New York: Springer, 1996

\bibitem {P69} Putnam H.  An unsolvable problem in number theory.
 J. Symbolic Logic, 1960, 25: 220--232

\bibitem {Ri89} Ribenboim P. The Book of Prime Number Records, 2nd Edition.
 New York: Springer, 1989

\bibitem {Ro49} Robinson J.  Definability and decision problem in arithmetic.
 J. Symbolic Logic, 1949, 14: 98--114

\bibitem{Sh} Shlapentokh A. Hilbert's Tenth Problem: Diophantine Classes and Extensions to Global Fields. New Mathematical Monographs, Vol. 7, Cambridge: Cambridge Univ. Press, 2007.

\bibitem{S-PhD} Sun Z-W. Further results on Hilbert's tenth problem (in Chinese).
PhD Thesis. Nanjing: Nanjing University, 1992

\bibitem {S92a} Sun Z-W.  Reduction of unknowns in Diophantine representations.
Sci. China Ser. A, 1992,  35: 257--269. Available from {\tt http://maths.nju.edu.cn/$\sim$zwsun/12d.pdf}

\bibitem {S92b} Sun Z-W.  A new relation-combining theorem and its application.
 Z. Math. Logik Grundlag. Math., 1992, 38: 209--212

\bibitem {S15} Sun Z-W. On universal sums of polygonal numbers. Sci. China Math.,
 2015, 58: 1367--1396

\bibitem {S16} Sun Z-W.  A result similar to Lagrange's theorem. J. Number Theory,
2016, 162: 190--211

\bibitem {T85} Tung S P,  On weak number theories. Japan. J. Math. (N.S.), 1985, 11: 203--232

\bibitem {T87} Tung S P,  Computational complexities of Diophantine equations with parameters. J. Algorithms, 1987, 8: 324--336

\bibitem {V97} Vaughan R C. The Hardy-Littlewood Method, 2nd Edition. Cambridge Tracts in Math., Vol. 125, Cambridge: Cambridge Univ. Press, 1997


\end{thebibliography}
 \end{document}